%% file: agt-3-36.tex
\documentclass{gtart}
\input agtout

\lognumber{36}
\volumenumber{3}
\volumeyear{2003}
\papernumber{36}
\published{16 October 2003}
\pagenumbers{1051}{1078}
\received{17 March 2003}
\revised{15 September 2003}
\accepted{21 September 2003}

\usepackage{amsfonts}
\usepackage{amsmath}
\usepackage{amssymb}
\usepackage{eufrak}
\usepackage{epsfig}
\usepackage[all]{xy}
\usepackage{longtable}
\usepackage{mathrsfs}

\newcommand{\lra}{\longrightarrow}
\newcommand{\Rem}{\par\medskip {\bf Remark}\qua}
\newcommand{\Cl}{\par {\bf Claim}\qua}
\newcommand{\Def}{\par\medskip {\bf Definition}\qua}

\let\Bbb\mathbb

\newcommand{\Par}{\partial}
\newcommand{\lbar}[1]{\overline{#1}}
\newcommand{\inter}[1]{\mathring{#1}}

\newcommand{\id}{\mbox{id}\;}
\renewcommand{\mod}{\mbox{mod}\;}
\newcommand{\sign}{\mbox{sign}}
\newcommand{\im}{\mbox{im}\;}
\newcommand{\coker}{\mbox{coker}}

\newcommand{\pr}{\text{pr}}
\renewcommand{\b}[1]{{\mathbf #1}}
\renewcommand{\phi}{\varphi}

\renewcommand{\cal}{\mathcal}

\newcommand{\f}{\mathscr}

\begin{document}
\title{Resolutions of p-stratifolds with\\isolated singularities}
\author{Anna Grinberg}
\address{Department of Mathematics, UC San Diego\\9500 
Gilman Drive, La Jolla, CA, 92093-0112}
\email{agrinber@math.ucsd.edu}

\begin{abstract}{Recently M. Kreck introduced a class of
stratified spaces called p-stratifolds \cite{Kr3}.  He defined and
investigated resolutions of p-stratifolds analogously to resolutions
of algebraic varieties. In this note we study a very special case of
resolutions, so called optimal resolutions, for p-stratifolds with
isolated singularities. We give necessary and sufficient conditions
for existence and analyze their classification.}
\end{abstract}

\asciiabstract{Recently M. Kreck introduced a class of
stratified spaces called p-stratifolds [M. Kreck, Stratifolds,
Preprint]. He defined and investigated resolutions of p-stratifolds
analogously to resolutions of algebraic varieties. In this note we
study a very special case of resolutions, so called optimal
resolutions, for p-stratifolds with isolated singularities. We give
necessary and sufficient conditions for existence and analyze their
classification.}

\maketitle

\newtheorem{thm}{Theorem}%[section]
\newtheorem{lem}[thm]{Lemma}
\newtheorem{cor}[thm]{Corollary}
\newtheorem{prop}[thm]{Proposition}
\newtheorem*{cori}{Corollary}

\let\\\par

\section{Introduction}\label{int}
Roughly speaking, p-stratifolds are  topological spaces which are
constructed by attaching manifolds with boundary by a map to the
already inductively constructed space. The attaching map has to
fulfill some subtle
properties. There is a more general notion of stratifolds introduced
by M. Kreck \cite{Kr3}. However, the only results concerning
the resolution of stratifolds exist after going over to the subclass 
of p-stratifolds.\\

The situation simplifies very much, if we consider only
p-stratifolds with isolated singularities, where the construction
is done in two steps only. The first step is the choice of a
countable number of points $\{x_i\}_{i\in I\subseteq \mathbb N}$
which will become the isolated singularities. The second step is
the choice of a smooth manifold $N$ of dimension $m$, together
with a proper map $g: \Par N \lra \{x_i\}_{i\in I}$, where $\{x_i\}_{i\in I}$ is
considered as
    $0$-dimensional manifold and the collection of boundary components
    $f^{-1}(x_i)$ is equipped with a germ of collars. The p-stratifold is obtained by forming
$$
\f S = N \cup_g \{x_i\}_{i\in I}.
$$
 We reformulate this in a slightly different way.

\Def An $m$-dimensional  {\em p-stratifold with isolated
singularities} is a topological space $\f S$ together with a
proper map $f: N \lra \f S$, where
\begin{itemize}
\item $N$ is an $m$-dimensional manifold with boundary,
\item $f|_{\inter N}$ is  a homeomorphism onto its image,
\item $\f S - f(\inter N)$ is a discrete countable set, denoted by $\Sigma$,
  the {\em singular set},
\item $f^{-1}(x)$ is equipped with a germ of
  collars for all $x \in \Sigma$,
\item $U \subset \f S$ is open if and only if $U \cap \Sigma$ is
open
  and $f^{-1}(U)$ is open in $N$.
\end{itemize}
The manifold $f(\inter N)$ is called the {\em top stratum},
$\Sigma$ is called the $0$-stratum of $\f S$. Choose an
identification $\Sigma = \{x_i\}_{i\in I\subseteq \Bbb N}$ and
denote the collection of boundary components mapped to a singular
point $L_i:=f^{-1}(x_i)$  the
{\em link of $\f S$ at $x_i \in \Sigma$}.\\

A {\em collar} around $L_i$ is a diffeomorphism 
$\b c_i: L_i \times [0,\epsilon_i) \lra U_i$, where $U_i$ is an open
neighbourhood of $L_i$ in $N$ and $\epsilon_i >0$, such that
$c_i|_{L_i \times \{0\}}$ is the identity map on $L_i$. The {\em germ}
is an equivalence class of collars, where two collars
$\b c_i: L_i \times [0,\epsilon_i)\lra U_i$ and 
$\tilde{\b c}_i: L_i \times [0,\tilde{\epsilon}_i)\lra \tilde{U}_i$
are called equivalent if there is a positive 
$\delta \leq \min\{\epsilon_i,\,\tilde{\epsilon}_i\}$, such that
$\b c_i|_{L_i\times [0,\delta)} \equiv \tilde{\b c}_i|_{L_i\times [0,\delta)}$.
The role of the collars becomes clear if we define smooth maps
from a smooth manifold to a p-stratifold.

\Def Let $\f S$ be a p-stratifold with isolated singularities  and
$\rho: \f S \lra \Bbb R$ a continuous map. The map $\rho$ is
called {\em smooth} if $\rho f|_{\inter N} :\inter N \lra \Bbb R$
is smooth and there are representatives of the germ of collars $\b
c_i: L_i \times [0,\varepsilon_i) \lra N$  satisfying for all $i$:
$$
\rho f(\b c_i(x,t)) = \rho f(x) \mbox{ for all } x \in L_i.
$$
Let $M$ be a smooth manifold. A continuous
 map $g: M \lra \f S$ is called {\em smooth}, if for all smooth
maps $\rho : \f S \lra \Bbb R$ the composition $\rho g$ is
again smooth.\\
\noindent
It is not hard to verify that the  map $g$ is smooth
if and only if the restriction
$g|_{g^{-1}(\f S -  \Sigma)}$ is smooth.\\

The most important examples of p-stratifolds as defined above are
algebraic varieties with isolated singularities.\\

\medskip
{\bf Example}\qua (Algebraic varieties with isolated singularities)\\
Consider an algebraic variety $V \subset \Bbb R^n$ with isolated
singularities, i.e.\ the singular set $\Sigma$ is zero-dimensional.
%$\{s_1,\dots,s_k\}$, i.e.\ $dp(s_i) = 0$, but there is an
%$\varepsilon > 0$ such that $dp(x) \neq 0$ for all $x \in
%(p^{-1}(0)\cap D_{\varepsilon}(s_i)) - \{s_i\}$, where
%$D_{\varepsilon}(s_i) := \{x \in \Bbb R^{n+1}|\; ||s_i -x|| \leq
%\epsilon\}$.\\
Let $s_i \in \Sigma$ be a singular point. There is nothing to do
if $s_i$ is open in $V$. Otherwise consider the distance function
$\rho_{i}$ on $\Bbb R^{n}$ given by  $\rho_i(x) := ||x-s_i||^2$.
It is well known that there is an $\varepsilon_i >0$ such that on
$V_{\varepsilon_i}(s_i) := V \cap D_{\varepsilon_i}(s_i)$ the
restriction $\rho_i|_{V_{\varepsilon_i}-\{s_i\}}$ has no critical
values. Here $D_{\varepsilon_i}(s_i)$ denotes the closed ball in
$\Bbb R^n$ of radius $\varepsilon_i$ centered at $s_i$. Set $\Par
V_{\varepsilon_i}(s_i) := V_{\varepsilon_i}(s_i) \cap \Par
D_{\varepsilon_i}(s_i)$. \\

By following the integral curves of the gradient vector field of
$\rho_i|_{V_{\varepsilon_i}-\{s_i\}}$, we obtain a diffeomorphism
$$
\xymatrix{ h: \Par V_{\varepsilon_i}(s_i)\times [0,\epsilon_i)
\ar[rr]\ar[rd]_{\pr_2}
& &V_{\varepsilon_i}-\{s_i\} \ar[ld]^{\epsilon_i - \rho_i}\\
& [0,\epsilon_i) }
$$
being the identity on $\Par V_{\varepsilon_i}(s_i)\times \{0\}$,
see \cite[\S 6.2]{H}. We extend this map to a continuous map
$$
\bar h: \Par V_{\varepsilon_i}(s_i) \times [0,\epsilon_i] \lra
V_{\varepsilon_i}.
$$
Finally, we define the manifold $N$ (with obvious collar) by
setting
$$
N := V - (\sqcup_i \inter{D}_{\varepsilon_i}(s_i))\cup_{\id} \Par
V_{\varepsilon_i}(s_i)\times [0,\epsilon_i].
$$
The map $f = \id\cup \bar h: N \lra V$ gives $V$ the structure of
a p-stratifold  with isolated singularities.

Since every complex algebraic variety is in particular a real one,
we obtain the same result for a complex algebraic variety with
isolated singularities.\\

From now on all p-stratifolds are p-stratifolds with isolated
singularities. To simplify the notation combine the
representatives of the collars $\b c_i: L_i
\times[0,\varepsilon_i)$ to a single map $\b c: \sqcup_i L_i
\times[0,\varepsilon_i) \lra N$.
%top stratum $f(\inter N)$
Using this map the singular set $\Sigma$ is equipped with the germ
of neighbourhoods $[U\Sigma]$ by taking $U\Sigma :=  f(\im\b
c)\sqcup (\Sigma-f(\Par N))$. The collars also give us a
retraction $r:U\Sigma \lra \Sigma$.

We also introduce the germ of closed neighbourhoods $[\lbar{U}
\Sigma]$ by setting $\lbar{U}\Sigma := f(\b c(\sqcup _i L_i \times
[0,\varepsilon_i/2])) \sqcup (\Sigma-f(\Par N))$. If we want to
make the dependency on the representative of the germ of collars
clear, we sometimes write $U\Sigma_{\b c}$ and $\lbar{U}\Sigma_{\b
c}$.
%By
%making this construction inductively, we obtain a germ of
%neighbourhood for each $k$-skeleton  $\Sigma_{k}$ ($k < m$). We denote
%the codimension 1 skeleton $\Sigma_{m-1}$ by $\Sigma$.\\

\Def Let $\f S$ be an $n$-dimensional p-stratifold with top
stratum $f(\inter N)$. A {\em resolution} of $\f S$ is a proper
map $p: \hat{\f S} \lra \f S$ such that
\begin{itemize}
\item $\hat{\f S}$ is a smooth manifold; \item $p$ is a proper
smooth map;
%\item $p$ is a proper smooth map;
\item the restriction of $p$ on $p^{-1}(f(\inter N))$ is a
  diffeomorphism on $f(\inter N)$;
\item $p^{-1}(f(\inter N))$ is dense in $\hat{\f S}$;
\item the
inclusion $ \hat\Sigma := p^{-1}(\Sigma) \hookrightarrow
  U\hat{\Sigma} := p^{-1}(U\Sigma)$ is a homotopy equivalence for a
  representative of the neighbourhood $U\Sigma$ of $\Sigma$.
\end{itemize}
A resolution $p: \hat{\f S} \lra \f S$ is called {\em optimal}, if
$p|_{\hat \Sigma}:\hat \Sigma \lra \Sigma$ is an
$[n/2]$-equivalence. In particular, it follows that $p: \hat{\f S}
\lra \f S$ is an
$[n/2]$-equivalence as well.\\

If the manifold $N$ is equipped with more structure, e.g.\
orientation  or spin-structure, we introduce corresponding
resolutions, which have more structure.

\Def Let $\f S = f(N) \cup \{x_i\}_i$ be a p-stratifold with
oriented $N$. A resolution $p: \hat{\f S} \lra \f S$ is called an
{\em oriented resolution}, if $\hat{\f S}$ is oriented and
$p|_{p^{-1}(f(\inter{N}))}$ is orientation preserving.
Analogously, if $N$ is spin, then $p: \hat{\f S} \lra \f S$ is
called a {\em spin resolution} if $\hat{\f S}$ is
spin and $p|_{p^{-1}(f(\inter{N}))}$ preserves the spin structure.\\

If $V$ is an algebraic variety, Hironaka has shown \cite{Hi} that
there is a resolution of singularities in the sense of algebraic
geometry. The above topological definition is modelled on the one
from algebraic geometry.  All conditions are analogous except the
last one, which is  always fulfilled in the context of algebraic
geometry. As explained in \cite{B-R}, a neighbourhood $U$ of the
singular set $\Sigma$ of an algebraic variety $V$ such that the
inclusion $\Sigma \hookrightarrow U$ is a homotopy equivalence can
be obtained from a proper algebraic map $\rho: V \lra \Bbb R$ with
$\Sigma = \rho^{-1}(0)$ by taking $U = \rho^{-1}[0,r)$, provided
$r>0$ is small enough. Thus for a resolution  $p: \hat V \lra V$
the preimage $\hat U := p^{-1}(U)$ is a neighbourhood of $\hat
\Sigma:=p^{-1}(\Sigma)$ in $\hat V$ obtained from $\hat \rho :=
\rho p$, hence the inclusion $\hat \Sigma
\hookrightarrow \hat U$ is a homotopy equivalence.\\

 Note that $\lbar{U}\hat\Sigma$ is a smooth manifold with boundary
diffeomorphic to $\Par N$. Consider  the preimage of the
neighbourhood of each singularity and set $\lbar{U}\hat\Sigma_i :=
p^{-1}(f{\b c}_i(L_i \times [0,\varepsilon_i/2]))$, where ${\b
c}_i$ is the restriction of the collar
to $L_i \times [0,\varepsilon_i/2]$.\\
It is not hard to verify that
 a resolution $p: \hat{\f S} \lra \f S$ is
optimal if and only if the manifolds $\lbar{U} \hat \Sigma_i$ are
$([n/2]-1)$-connected.\\

%After Hironaka's great result a resolution of an algebraic
%variety  always exists, whereas  we obtain obstructions for the existence
%of resolutions.
In contrast to algebraic varieties, resolutions of stratifolds in
general do not exist, not even for isolated singularities. But in
this case there is a simple necessary and sufficient condition,
see \cite{Kr3} and \S \ref{proofs} for a proof.

\begin{thm}
\label{k} An $n$-dimensional p-stratifold with isolated
singularities
% ${x_1,\dots,x_k}$. Then $S$
admits a resolution if and only if each link of the singularity
$L_i$,
%for $i\in {1,\dots,k}$
vanishes
in the bordism group $\Omega_{n-1}$.
\end{thm}

\noindent {\bf Example}\qua The p-stratifold $\f S= \Bbb CP^2 \times
I \cup_f \{x_0,x_1\}$ with the obvious stratification, such that
$f(\Bbb CP^2 \times \{0\}) = x_0$ and $f(\Bbb CP^2 \times \{1\}) =
x_1$, does not admit a resolution.\\

%Even in the case of isolated singularities, there are too many
%resolutions, hence the classification question is very hard to
%answer. We shrink our considerations to a special class of
%resolutions.\\
To give a feeling of the result concerning optimal resolution, we
formulate the following  special case which will be derived as
Corollary \ref{7} of Theorem \ref{g2} (cf.\ \S \ref{P2}).

\begin{cori}
Let $\f S$ be a p-stratifold with parallelizable links of
singularities $L_i$. Assume $L_i$ is bounded by a parallelizable
manifold, then $\f S$ admits an optimal resolution.
\end{cori}

We have shown above that every algebraic variety with isolated
singularities admits a structure of a p-stratifold. One may ask
the  converse question. When does a p-stratifold with isolated
singularities admit an algebraic structure? The following Theorem
of Akbulut and King \cite[Thm. 4.1]{AK} clarifies the situation in
the case of a real algebraic structure.
\begin{thm}
A topological space $X$ is homeomorphic to a real algebraic set
with isolated singularities if and only if $X$ is obtained by
taking a smooth compact manifold $M$ with boundary $\Par M =
\cup_{i=1}^r L_i$, where each $L_i$ bounds, then crushing some
$L_i$'s to points and deleting the remaining $L_i$'s.
\end{thm}
Combining this result with Theorem \ref{k} we immediately obtain:
\begin{cor}
A compact p-stratifold $\f S$ with isolated singularities is
homeomorphic to a real algebraic set with isolated singularities
if and only if $\f S$ admits a resolution.
\end{cor}

\noindent {\bf Example}\qua (Resolutions of hypersurfaces with
isolated
  singularities)

Let $p: \Bbb R^{n+1} \lra \Bbb R$ be a polynomial with isolated
singularities $\{s_i\}_i$, i.e.\ $s_i \in V := p^{-1}(0)$ and $s_i$
is an isolated critical point of $p$.  Assume further that the
points $s_i$ are not open. According to a previous  example, the
hypersurface $V$ admits a canonical structure of a p-stratifold.
We have to investigate the link of the singularity, which is given
by $\Par V_{\varepsilon_i}(s_i)$.

Choose a $\delta > 0$ such that all $c$ with $|c| \leq \delta$ are
regular values of $p$ and take $c$ such that $p^{-1}(c) \neq
\emptyset$.  Then $p^{-1}(c)$ is a smooth manifold with trivial
normal bundle. With the help of the gradient vector field we see
that $p^{-1}(c)\cap S^{n}_{\varepsilon_i}(s_i)$ is diffeomorphic
to $p^{-1}(0) \cap S^{n}_{\varepsilon_i}(s_i) = \Par
V_{\varepsilon_i}(s_i) $. Thus, $\Par V_{\varepsilon_i}(s_i)\cong
p^{-1}(c)\cap S^{n}_{\varepsilon_i}(s_i) \cong   \Par
(p^{-1}(c)\cap D^{n+1}_{\varepsilon_i}(s_i))$.
%is a boundary of a parallelizable manifold.
We see that a resolution always exists, and since the bounding
manifolds are automatically parallelizable, we  even obtain an
optimal resolution after choosing an appropriate bordism (compare
with Figure \ref{fig1}).
\\
%In the real case $\Par V_{\varepsilon}$ is a codimension 1 oriented
%submanifold of $S_{\varepsilon}$, thus the Jordan curve theorem implies that
%all components of $\Par V_{\varepsilon}$ bound in $S_{\varepsilon}$. Hence a
%resolution always exists, and since the bounding manifolds are automatically stabile
%parallelizable, we  even obtain an optimal resolution.
\setlength{\unitlength}{1cm}
\begin{figure}[ht!]\small
\begin{picture}(12,6)
\put(3,0){\epsfig{file=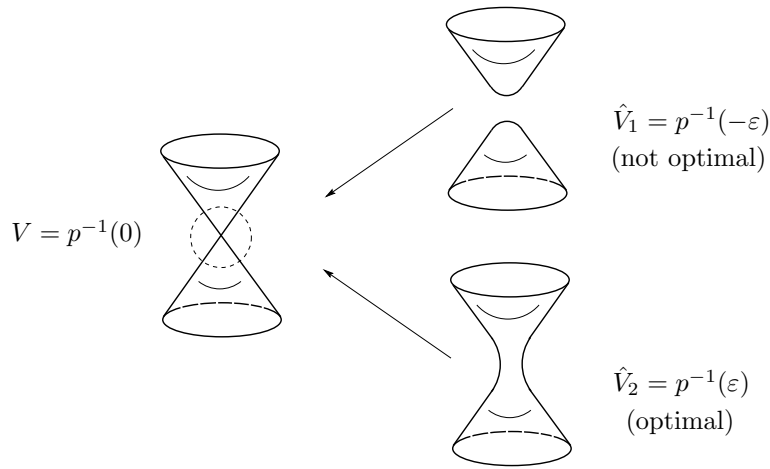,width=5.5cm}}
%\put(0,5.5){$p(x,y,z) = x^2 +y^2 - z^2$}
\put(1,3){$V = p^{-1}(0)$} \put(9,4.5){$\hat V_1 = p^{-1}(-
\varepsilon)$} \put(9,4){(not optimal)} \put(9,1){$\hat V_2 =
p^{-1}(\varepsilon)$} \put(9.2,0.5){(optimal)}
\end{picture}
\caption{\label{fig1}$p(x,y,z) = x^2 +y^2 - z^2$}
\end{figure}

In the case of a complex polynomial $p: \Bbb C^{n+1} \lra \Bbb C$
$(n>0)$, every deformation $p^{-1}(c)$ gives us an optimal
resolution, provided $|| c||$ is small enough.
%To see this note first since zero is the only critical value
%of $p$ near zero  there is an $\epsilon > 0$ such that for
%$||c|| < \epsilon$ the preimage $p^{-1}(c)$ is smooth and has trivial
%normal bundle. Further with help of a gradient vector field we see
%that   $p^{-1}(c)\cap S^{2n+1}_{\varepsilon}$ is diffeomorphic to
%$p^{-1}(0) \cap S^{2n+1}_{\varepsilon}$. Thus
%$L:= p^{-1}(0)\cap S^{2n+1}_{\varepsilon} \cong p^{-1}(c)\cap
%S^{2n+1}_{\varepsilon} \cong   \Par (p^{-1}(c)\cap
%D^{2n+2}_{\varepsilon})$ parallelizable. Now we
This follows from a result of Milnor \cite[Thm. 6.5]{Mi} which
states that $M := p^{-1}(c)\cap D^{2n+2}_{\varepsilon}(s_i)$ is
homotopy equivalent to a wedge of $\mu \geq 1$ copies of $S^n$ and
thus $(n-1)$-connected.\\
%------------------------------------------------------------------------

Consider another interesting class of p-stratifolds with isolated
singularities, namely those arising from a smooth group action.

\Def A smooth $S^1$-action on a smooth manifold $M$ is called {\em
semi-free} if the action is free outside of the fixed point set,
i.e.\ if $gx  = x$ for a $g \in S^1, g \neq 1$ and $x \in M$, then
$hx = x$ for all $h\in S^1$.

\begin{lem}
Let $M$ be a closed oriented manifold with semi-free $S^1$-action
with only isolated fixed points. Then $M/S^1$ admits an oriented
resolution if and only if $\dim M \equiv 0 (\mod 4)$.
\end{lem}
\begin{proof}

Let $\dim M = n$. There is nothing to show if the action is free.
Thus let $x\in M$ be a fixed point. The differential of the action
gives a representation of $S^1$ on $T_xM$ and there is an
equivariant local diffeomorphism from $T_xM$ onto a neighbourhood
of $x$ in $M$. According to \cite[Prop. (II.8.1)]{BT}, every
irreducible representation of $S^1$ on $\Bbb R^n$ is equivalent
to:
$$
\begin{array}{ccc}
{\begin{pmatrix}
z^{n_1}&  0 & \dots & \dots &  0 \\
0 & \cdot & \ddots&  & \vdots\\
\vdots&\ddots  & \cdot &\ddots  & \vdots\\
\vdots & &\ddots  & z^{n_{[n/2]}} & 0\\
0 &  \dots& \dots & 0& 1
\end{pmatrix}} &  &
%--------------------------------------------
{\begin{pmatrix}
z^{n_1}&  0 & \dots & \dots &  0 \\
0 & \cdot & \ddots&  & \vdots\\
\vdots&\ddots  & \cdot &\ddots  & \vdots\\
\vdots & &\ddots  & \cdot & 0\\
0 &  \dots& \dots & 0&  z^{n_{n/2}}
\end{pmatrix}}
\\
\\
\mbox{ considered as a representation on}&  & \mbox{ considered as a representation on} \\
\Bbb C \times \dots \times \Bbb
  C \times \Bbb R, \mbox{if $n$ is odd} & & \Bbb C \times \dots \times \Bbb
  C , \mbox{if $n$ is even}\\
\end{array}
$$

Since the action is semi-free and $x$ an isolated fixed point we
conclude that $\dim M$ is even. We can further assume $n_i=1$ for
all $i \in \{1,\dots,n/2\}$. Let $\dim M = 2m$ and let
$\{x_1,\dots,x_k\}$ be the set of fixed points. Choose equivariant
disks $D_{x_i}$ around $x_i$. In this situation we have
$$
M/S^1 = (M - \sqcup_i \inter{D}_{x_i})/S^1 \cup \{x_1,\dots,x_k\}.
$$
The domain of  the top stratum is then given by $ N := (M -
\sqcup_i \inter{D}_{x_i})/S^1$ and the singular set is $\Sigma :=
\{x_1,\dots,x_k\}$. The links of singularities are given by $L_i
\cong S^{2m-1}/S^1 = \Bbb C P^{m-1}$. Using Theorem \ref{k} we
conclude that the resolution exists if and only if $[\Bbb
CP^{m-1}]$ vanishes in $\Omega^{SO}_{2m-2}$. For $m = 2l+1$ the
signature of $\Bbb CP^{m-1}$ is equal to 1, hence $\Bbb CP^{m-1}$
does not bound. In the case of an even $m = 2l+2$ we have $\Bbb
CP^{2l+1} = S^{4l +3} / S^1 = S(\Bbb H^{l+1})/S^1$, where $S(\Bbb
H^{l+1})/S^1$ is the sphere bundle
$$
\xymatrix{
S^2 = S^3/S^1 \ar@{^(->}[r] & S(\Bbb H^{l+1})/S^1 \ar[d] \\
& S(\Bbb H^{l+1})/S^3 = \Bbb HP^l }
$$
and the associated disk bundle bounds.
\end{proof}

%----------------------------------------------------------------------------

%\section{Equivalent Resolutions}
As mentioned before, we are particularly interested in the
classification of resolutions. Thus we have to decide when we are
going to consider two resolutions as equivalent. We can restrict
our attention to the resolving manifolds and introduce a relation
on them, e.g.\ diffeomorphism, but in this case we completely
ignore an important part of the resolution data, namely the
resolving map. Hence, one can ask for diffeomorphisms between the
resolving manifolds commuting with the resolving maps. This
relation is  very strong and, therefore, very hard to control. In
the following
definition, we combine these two ideas.\\

\Def Let $\f S$ be a p-stratifold and $p: \hat{\f S} \lra  \f S$
and $p': \hat{\f S}' \lra \f S$ two resolutions of $\f S$. We call
the resolutions {\em
  equivalent}, if, for  every representative of the neighbourhood germ
$\lbar{U}\Sigma_\b c$, there is a diffeomorphism $\phi_{\b c}:
\hat{\f S} \lra \hat{\f S}'$ such that  the following holds:
\begin{itemize}
\item $p'\phi_\b c = p$ on $\hat{\f S} - \lbar{U}\hat\Sigma_\b c$
and \item $r p'\phi_\b c = rp$ on $\lbar{U}\hat\Sigma_\b c$, where
$r:\lbar{U}\Sigma_\b c \lra
  \Sigma$ is the neighbourhood's retraction.
\end{itemize}

This means outside of an arbitrary small neighbourhood of the
singularity, the diffeomorphism commutes with the resolving maps
and near $\Sigma$ it only commutes after the composition with the
retraction.

Observe that $\phi_\b c$ gives a diffeomorphism $\Par \lbar{U}\hat
\Sigma_{\b c} \lra \Par \lbar{U} \hat \Sigma'_{\b c}$.\\

The classification of optimal resolutions is quite a difficult
problem. For if $\hat{\f S}$ is an optimal resolution of $\f S$,
then $\hat{\f S} \sharp \mathcal S$ is again optimal for an
arbitrary homotopy sphere $\mathcal S$. In particular, consider
the sphere $S^n$ stratified as $D^n \cup {pt}$, then every
homotopy sphere $\mathcal S^n$ gives us a resolution of $S^n$.
Thus, we weaken the problem and ask for the equivalence up to a
homotopy sphere.\\

\Def Two resolutions  $\hat{\f S} \lra  \f S$ and $\hat{ \f S}'
\lra \f S$ are called {\em almost equivalent} if $\hat{\f S}
\sharp \mathcal S$
is equivalent to $\hat{\f S}'$ for a homotopy sphere $\mathcal S$. \\

A special case of the classification result is the following.
\begin{cori}
Let $\hat{\f S} \lra \f S$ and $\hat{\f S'} \lra \f S$ be two
resolutions of a $2n$-dimens\-ion\-al p-stratifold $\f S$
%$\{x_1,\dots,x_l\}$
having $(n-2)$-connected links of isolated singularities. Assume
that $n \equiv 6\; (\mod 8)$and that $\lbar{U} \hat \Sigma_i$ and
$\lbar{U} \hat \Sigma'_i$ are parallelizable with compatible
parallelizations on the boundary. Let further $e(\lbar{U} \hat
\Sigma_i)= e(\lbar{U} \hat \Sigma'_i)$ and $\sign(\lbar{U} \hat
\Sigma_i \cup_\Par \lbar{U} \hat \Sigma'_i) =0$. Then there is a
$k\in \{0,1\}$ such that $\hat{\f S} \sharp k(S^{n} \times S^{n})$
is almost equivalent to $\hat{\f S}'\sharp k(S^{n} \times S^{n})$.
\end{cori}

\section{Existence of optimal resolutions}\label{P2}

Before proceeding with the existence of an optimal resolution we
need to introduce some notation. For a topological space $X$ let
$X \langle k \rangle$ be the $k$-connected cover of $X$, which
always comes with a fibration $p:X \langle k \rangle \lra X$. For
further information see for example \cite{Ba}. We take $X$ to be
the classifying space $BO$ and denote by $\Omega_n^{BO\langle k
\rangle}$ the bordism group of closed $n$-dimensional manifolds
together with a lift of the normal Gauss map, compare \cite[Chap.
I]{St}.

\begin{thm} \label{g2}
An $n$-dimensional  p-stratifold with isolated singularities
admits an optimal resolution if and only if the normal Gauss map
$\nu_j : L_j \lra BO$ admits a lift over
$BO\langle[n/2]-1\rangle$, such that $[L_j,\bar{\nu}_j] =0$ in
\mbox{$\Omega_{n-1}^{BO\langle[n/2]-1\rangle }$}.
\end{thm}

\begin{proof}
Let $\f S$ be a p-stratifold with isolated singularities
$\{x_i\}_{i\in I}$ and $\hat{\f S}$ an optimal resolution of $\f
S$. Set $W_i := \lbar{U}\hat\Sigma_i$, then  $W_i$ is a smooth
manifold with boundary $L_i$ for $i \in I \subseteq \Bbb N$.
Consider the normal Gauss map, together with the
$([n/2]-1)$-connected cover over $BO$.
$$
\xymatrix{
 &
 BO\langle [n/2]-1\rangle  \ar[d] &  F \ar@{_{(}->}[l] \\
{ \nu_j : {W}_j \ar[r] \ar@{-->}[ur]^{\bar{\nu}_j} } &  BO
}
$$

The obstructions for the existence of a lift lie in
$\tilde{H}^r({W}_j,\pi_{r-1}(F))$. Note that we can use global
coefficients since the fibration  $BO\langle [n/2]-1\rangle \lra
BO$ is
simple.\\

Since the resolution is optimal the manifold $W_i$ is
$([n/2]-1)$-connected, hence $\tilde{H}^r({W}_i)=0$ for $r <
[n/2]$.

Using the properties of the connected cover it follows from the
long exact homotopy sequence %for $BO\langle [n/2]-1\rangle \lra BO$
%together with the properties of the connected cover it follows
that $\pi_r(F) = 0$ for $r \geq [n/2]-1$.\\

Hence, there are no obstructions for the lifting of the normal
Gauss map, thus $[L_j,\nu_j]$ vanishes in
$\Omega_{n-1}^{BO\langle[n/2]-1 \rangle }$.

The fact that the condition is also sufficient is an immediate
consequence of the following result from \cite[Prop. 4]{Kr1}.
\end{proof}

\begin{thm}\label{k2}
Let $\xi : B \lra BO$ be a fibration and assume that $B$ is connected
and has a finite $[m/2]$-skeleton. Let $\bar \nu : M \lra B$ be a
lift of the normal Gauss map of an $m$-dimensional compact manifold $M$. Then
if $m \geq 4 $, by a finite sequence of surgeries $(M, \bar \nu)$ can
be replaced by $(M',\bar{\nu}')$ so that $\bar{\nu}': M' \lra B$ is an
$[m/2]$-equivalence.
\end{thm}

For example, we obtain the following:
\begin{cor}\label{7}
Let $\f S$ be a p-stratifold with parallelizable links of
singularities $L_i$. Assume $L_i$ is bounded by a parallelizable
manifold, then $\f S$ admits an optimal resolution.
\end{cor}

%\begin{lem}
%Let $\Phi$ denote the field of real or complex numbers and let $n>1$. Each
%hypersurface $V = p^{-1}(0)$ in $\Phi^{n+1}$ with isolated
%singularities admits an optimal resolution. In the case $\Phi = \Bbb
%C$, the deformation $p^{-1}(c)$ is itself the total space of an
%optimal resolution, provided $||c||$ is small enough.
%\end{lem}

%---------------------------------------------------------------------------

%\Rem Let $p : \hat S \lra S$ be a resolution of a p-stratifold
%$S$. We can write $S$ as $S_n \cup_{f|_{\Par S_n}}\Sigma_{n-1}(S)$, and
%according to the definition of resolution
%$\hat S \cong S_n \cup_{\hat f|_{\Par S_n}}\Sigma_{n-1}(\hat S)$.
%Using  the collar of $S_n$ we see, that $\hat S$ is determined by
%$\hat W := \Par S_n \cup_{\hat f |_{\Par S_n}}\Sigma_{n-1}(\hat
%S)$. The space $W$ is a smooth manifold with boundary $\Par S_n$.\\

%
%identity on the boundary $\Par S_n$ can be
%extend to a diffeomorphism $\varphi$:
%$$
%\xymatrix{ \varphi : W_1 \ar[rr] \ar[dr]& &  W_2 \ar[dl] \\
% & S }
%$$

\section{Classification of almost equivalent resolutions}
\label{cl} Now we turn to the main result of this note. In this
section we consider p-stratifolds with isolated singularities of
dimension $2n > 4$ having $(n-2)$-connected links of
singularities. The  classification is based on the following
result from \cite[Thm. 2]{Kr1}.

\begin{thm}\label{k1}
For  $n > 2$ let $W_1$ and $W_2$ be two compact connected $2n$-manifolds with
normal $(n-1)$-smoothings in a fibration $B$. Let  $g: \Par W_1 \lra
\Par W_2$ be a diffeomorphism compatible with the normal
$(n-1)$-smoothings $\nu_1$ and $\nu_2$.
Let further
\begin{itemize}
\item $e(W_1) = e(W_2)$,
\item $[W_1 \cup_g (-W_2), \bar{\nu}_1\cup \bar{\nu}_2] = 0
  \in \Omega_{2n}(B)$.
\end{itemize}
 Then $g$ can be extended to a diffeomorphism
$G:W_1\sharp k(S^n \times S^n  ) \lra W_2\sharp k(S^n \times S^n
)$ for $k \in \Bbb N$. Moreover, if $B$ is 1-connected then $k\in\{0,1\}$.\\

If $W_1$ is simply connected and $n$ is odd, then $k$ can be chosen as
$0$, i.e.\ we obtain diffeomorphism instead of stable diffeomorphism.\\
\end{thm}

We have to explain some terms appearing in the last theorem. Let
$B$ be a fibration over $BO$, a {\em normal $B$-structure} on a
manifold $M$ is a lift $\bar \nu$ of the normal Gauss map $\nu : M
\lra
BO$ to $B$.\\

\Def Let $B$ be a fibration over $BO$.
\begin{enumerate}
\item A normal $B$-structure $\bar \nu: M \lra B$ of a manifold $M$ in
  $B$ is a {\em normal $k$-smoothing}, if it is a $(k+1)$-equivalence.
\item We say that $B$ is {\em $k$-universal} if the fiber of the map $B \lra
  BO$ is connected and its homotopy groups vanish in dimension $\geq k+1$.
\end{enumerate}
Obstruction theory implies that if $B$ and $B'$ are both $k$-universal
and admit a normal $k$-smoothing of the same manifold $M$, then the
two fibrations are fiber homotopy equivalent. Furthermore, the theory
of Moore-Postnikov decompositions implies that for each manifold $M$
there is a $k$-universal fibration $B^k$ over $BO$ admitting a normal
$k$-smoothing, compare \cite[\S 5.2]{Ba}. Thus, the fiber homotopy type of the
fibration $B^k$ over $BO$ is an invariant of the manifold $M$ and we
call it the {\em  normal $k$-type of $M$}.

There is an obvious bordism relation on closed $n$-dimensional
manifolds with normal $B$-structures and the corresponding bordism
group is denoted $\Omega_n(B)$.\\

Applying the theorem to our situation, we first  have to determine the
normal $(n-1)$-type of an $(n-1)$-connected $2n$ manifold.\\

Consider a subgroup $H$ of $G := \pi_n(BO)$. Since the last group is always
cyclic, the group $H$ is determined by an integer $k$, such that
$H=\langle kx \rangle$ where $x$ is the generator of $\pi_n(BO)$. We call this
integer the {\em index} of $H$.

Every subgroup $H$ of $G = \pi_n(BO)$ gives us a fibration
$$
\xymatrix{
B \ar[r] \ar[d]^{p} & P(K(G/H,n)) \ar[d] \\
BO\langle n-1\rangle  \ar[r]_{\theta} & K(G/H,n)
}
$$
where the map $\theta$ corresponds to the canonical epimorphism
$G \lra G/H$. We denote the space $B$ belonging to the
index-$k$ group $B_k$. The composition $p_k: B_k \stackrel{p}{\lra}
BO\langle n-1\rangle \lra BO$ gives us a fibration over $BO$ with
fiber $F_k$.\\

\Def A $2n$-dimensional manifold $M$ is said to have the {\em
index $k$}, if $\nu_*(\pi_n(M))$
is a subgroup of  index $k$ in $\pi_n(BO)$.\\

\begin{thm}\label{nt}
The fibration $B_k \lra BO$ is the normal $(n-1)$-type of an
$(n-1)$-connected $2n$-dimensional manifold $M$, if and only if
$M$ is of index $k$.
\end{thm}

%Apart from $(n-1)$-type and the Euler characteristic the bordism class
%of the glued manifold
%has to be considered,
%according to Theorem \ref{k1},
%this last condition  is hard to control. For example,
%if $\hat S$ is an optimal resolution of $S$, then $\hat S \sharp
%\mathcal S$ is again optimal for a homotopy sphere $\mathcal S$. In
%particular, consider the sphere $S^n$ stratified as $D^n \cup {pt}$,
%then every homology sphere $\mathcal S^n$ gives us a resolution of $S^n$.
%Thus we weaken the problem, and ask for the equivalence modulo
%homotopy sphere.\\

%\Def Tho resolutions $\hat S_1$ and $\hat S_1$ are called
%{\em quasi-equivalent}, if there is a homotopy sphere $\Sigma$, such that
%$\hat S_1\sharp\Sigma$ is equivalent to $\hat S_2$.\\
The proof can be found in \S \ref{proofs}, which also contains proofs of
the following two theorems.\\

Now we look for conditions implying an $(n-1)$-connected $2n$-manifold
to be bordant to a homotopy sphere. Note first that as an easy consequence
from the universal coefficient theorem the first
non-trivial homology
group is free. % and of even rank.
%these are simple  consequences
%from the universal coefficient theorem  and because of the intersection product of even rang.
The homological information of $M$ is stored in the triple
$(H_n(M),\Lambda , \nu)$ where $\Lambda$ denotes
the intersection product $ {\Lambda} : H_n(M) \lra \Bbb Z$,
we often
simply write $x\cdot y$ for $\Lambda(x,y)$.
The last data $\nu$ is the normal bundle information, described in the
following way.
According to a theorem of
Haefliger \cite{Hae} every element of $H_n(M)$ is represented by an
embedding $S^n \hookrightarrow M$, and two embeddings corresponding to the
same homotopy class are regular homotopic. Thus, assigning to an
embedded sphere its normal bundle gives us a well defined map
$\nu: H_n(M) \lra \pi_{n-1}(SO(n))$.\\

\Def
An $(n-1)$-connected $2n$-dimensional manifold $M$ is called
{\em elementary} if $H_n(M)$ admits a Lagrangian $L$ w.r.t.\ $\Lambda$, such that $\nu|_L
\equiv 0$.

\begin{thm}\label{hs}
Let $M$ be an $(n-1)$-connected manifold of dimension $2n$. Then $M$ is bordant
 to a homotopy sphere if and only if $M$ is elementary.
\end{thm}

%We call an $(n-1)$-connected $2n$-manifold $M$ fulfilling the
%conditions of the last Theorem {\em a manifold with Surgery-nice
%  homology}.\\

Our main result, based on the last two theorems is:

\begin{thm}\label{mt}
For $n>2$ let $\hat{\f S} \lra \f S$ and
  $\hat{\f S'} \lra \f S$ be two
optimal resolutions of a $2n$-dimensional p-stratifold $\f S$ with
isolated singularities $\{x_i\}_{i\in I}$, such that each link
$L_i$ is $(n-2)$-connected. Assume further that for a suitable
representative $\lbar{U}\Sigma$ the following conditions hold for
all $i \in I$:
\begin{itemize}
\item $e(\lbar{U}\hat \Sigma_i) = e(\lbar{U}\hat \Sigma'_i)$;
\item $\lbar{U}\hat \Sigma_i$ and $\lbar{U}\hat \Sigma'_i$ have the same index $k_i$;
\item there exits normal $(n-1)$-smoothings
$\bar \nu_i $ and ${\bar \nu'_i}$ of $\lbar{U}\hat \Sigma_i$ and $\lbar{U}\hat \Sigma'_i$ in the fibration $B_{k_i} \lra BO$, such
  that $\nu_i|_{\Par \lbar{U}\hat \Sigma_i} = \nu'_i|_{\Par \lbar{U}\hat \Sigma'_i }$;
\item $\lbar{U}\hat \Sigma_i \cup_{\Par}\lbar{U}\hat \Sigma'_i$ is elementary.
\end{itemize}
If $n$ is odd, then $\hat{\f S}$ is almost equivalent to $\hat{\f S'}$.\\
If $n$ is even, then $\hat{\f S}\sharp k(S^n\times S^n)$ is almost
equivalent to $\hat{\f S'}\sharp k(S^n\times S^n) $ for a $k \in
\{0,1\}$.
\end{thm}

\section{Algebraic invariants}
In this section we will find algebraic invariants, which allow us
to decide whether an $(n-1)$-connected closed $2n$-dimensional
manifold is
elementary or not $(n>2)$. Some proofs can be found in \S \ref{proofs_alg}. \\

Recall the algebraic data corresponding to such a manifold $M$. We
have a triple $(H,\Lambda,\nu_*)$, where $H=H_n(M)$ is a free
$\Bbb Z$-module, $\Lambda:H\times H \lra \Bbb Z$ is the
intersection product
%$(-1)^n$-symmetric unimodular
%quadratic form
and $\nu_*: H \lra \pi_{n-1}(SO_n)$ is a normal bundle map,
described in the previous section. The map $\nu_*$ is not a
homomorphism, but satisfies the following equation:
\begin{equation}
\nu_*(x+y) = \nu_*(x) + \nu_*(y) + \Par \Lambda(x,y), \tag{$*$}
\end{equation}
where $\Par :\Bbb Z \cong \pi_n(S^n) \lra \pi_{n-1}(SO_n)$ is the
  boundary map from the long exact homotopy sequence of the fibration
$SO(n) \hookrightarrow SO(n+1) \lra S^n$, see \cite{W0}.

Thus, we obtain an algebraic object, the set $\mathcal T_n$ of
triples $(H,\Lambda, \nu_*)$, where $H$ is a free $\Bbb Z$-module,
$\Lambda:H\times H \lra \Bbb Z$ is an $(-1)^n$-symmetric
unimodular quadratic form and $\nu_*: H \lra \pi_{n-1}(SO_n)$ is a
map satisfying $(*)$. We want to investigate the assumptions under
which  an element $(H,\Lambda,\nu_*)\in\mathcal T_n$ is
elementary, i. e. when  $H$ possesses a Lagrangian ${\cal L}$ with
respect
to $\Lambda$ such that $\nu_*|_{\cal L} \equiv 0$.  \\

We begin with an observation that for a $4k$-dimensional manifold,
the normal bundle information can be replaced by the stable
normal bundle map.
\begin{lem}
Let $n$ be even and let $S^n \hookrightarrow M^{2n}$ be an
embedding. The normal bundle $\nu(S^n)$ of $S^n$ in $M$ is trivial
if and only if $\nu \oplus \Bbb R$ is trivial and the Euler class
of $\nu(S^n)$ vanishes.
\end{lem}
Thus, instead of considering $(H,\Lambda,\nu_*)\in \mathcal T_n$
we can go over to $(H,\Lambda,s\nu_*)$, where $s\nu_*: H \lra
\pi_{n-1}(SO)$ corresponds to the stable normal bundle map. Since
the Euler class of an embedded sphere representing $x \in H$ can
be identified with the self intersection class we conclude:

\begin{lem} Let $n$ be even. Then
$(H,\Lambda,\nu_*) \in \mathcal T_n $ is elementary if and only if
$(H,\Lambda,s\nu_*)$ is elementary.
\end{lem}

Let $\mathcal T^s_n$ denote the set of triples
$(H,\Lambda,s\nu_*)$, with $H$ and $\Lambda$ as above and $s\nu_*:
H \lra \pi_{n-1}(SO)$ a homomorphism. According to the different
possibilities for
$\pi_{n-1}(SO)$ we distinguish 3 cases.\\

\noindent (1)\qua $\pi_{n-1}(SO)= 0$. \Cl $(H,\Lambda,s\nu_*) \in
\mathcal T^s_n$ is elementary if and only if $\sign(\Lambda)=0$,
where $\sign$ denotes
the signature of a quadratic form. \\

\noindent (2)\qua  $ \pi_{n-1}(SO) = \Bbb Z$. Since $\Lambda$ is
unimodular it induces an isomorphism $H \stackrel{\cong}{\lra}
H^*$, which we also denote by $\Lambda$. The map $s\nu_*$ gives an
element of $H^*$ and we consider $\kappa_{s\nu_*}:=
\Lambda^{-1}(s\nu_*) \in H$. \Cl $(H,\Lambda,s\nu_*) \in \mathcal
T^s_n$ is elementary if and only if
$\sign(\Lambda)=0$ and $\kappa_{s\nu_*}^2 = 0$.\\

\noindent (3)\qua  $\pi_{n-1}(SO) = \Bbb Z_2$. Let $(H,\Lambda,s\nu_*)$
be an element of $\mathcal T^s_n$ with vanishing signature and
suppose $\Lambda$ is of type $II$, i.e.\ $\Lambda(x,x) = 0 \;(\mod
2)$ for all $x \in H$. Note that since $n \neq 8$ in this case, an
elementary element corresponding to a manifold always has a type
$II$ quadratic form. Thus, the dimension of $H$ is even and
according to \cite[Lem. 9]{Mi_s} we can choose a basis
$\{\lambda_1,\dots,\lambda_k,\mu_1,\dots,\mu_k\}$ satisfying
$$
\Lambda(\lambda_i,\lambda_j) = 0, \;\;\; \Lambda(\mu_i,\mu_j)=0
\text{
  and }
\Lambda(\lambda_i,\mu_j) = \delta_{ij}.
$$
Consider the set of all elements $x \in H$ with $\Lambda(x,x)=0$
and denote its image under canonical projection on $H\otimes\Bbb
Z_2$ by $H^0$. The class $\Phi(H,\Lambda,s\nu_*) := \sum_{i=1}^k
s\nu_*(\lambda_i)s\nu_*(\mu_i)\in \Bbb Z_2$ is well-defined and is
equal to the value $s\nu_*$ takes most frequently on the finite
set $H^0$, the class is called Arf invariant.

\Cl An element $(H,\Lambda,s\nu_*) \in \mathcal T^s_n$ with type
$II$ form $\Lambda$ is elementary if and only if
$\sign(\Lambda)=0$ and
$\Phi(H,\Lambda,s\nu_*) = 0$.\\

Consider now the case of an odd $n$. The quadratic form is now
skew symmetric. Depending on the values of $\nu_*$ there are again
three
different cases (compare \cite{Ke}), which were completely investigated in \cite{W0}.\\

\noindent (4)\qua $\pi_{n-1}(SO_n) = 0$. In this case every element of
$\mathcal T_n$
is elementary.\\

\noindent (5)\qua $\pi_{n-1}(SO_n) = \Bbb Z_2$. As in $(3)$, we can
define the Arf invariant $\Phi(H,\Lambda,\nu_*) =
\sum_{i=1}^k\nu_*(\lambda_i)\nu_*(\mu_i) \in \Bbb Z_2$, using a
symplectic basis $\{\lambda_1,\dots,\lambda_k,\mu_1,\dots,\mu_k\}$
of $H$.

\Cl An element $(H,\Lambda,\nu_*) \in \mathcal T_n$ is elementary
if and
only if $\Phi (H,\Lambda,\nu_*)  = 0$.\\

\noindent (6)\qua $\pi_{n-1}(SO_n) = \Bbb Z_2\oplus \Bbb Z_2$. We
consider again the stable normal bundle map $s\nu_*: H \lra \Bbb
Z_2$, the projection on the first component. As in $(2)$ using
$\Lambda$, we obtain
 an element $\kappa$ (determined $\mod 2H$) with $s\nu_*(x) =
\Lambda(\kappa,x) \;(\mod 2)$ for all $x\in H$.

\Cl An element $(H,\Lambda,\nu_*)\in \mathcal T_n$ is elementary
if and only if $\Phi(H,\Lambda,\nu_*)=0$ and
$\mbox{pr}_2\nu_*(\kappa) = 0$,
where $\mbox{pr}_2$ denotes the projection on the second component.\\

Knowing the algebraic description of elementary manifolds, we
formulate a special case of Theorem \ref{mt}.

\begin{cor}
Let $\hat{\f S} \lra \f S$ and $\hat{\f S'} \lra \f S$ be two
resolutions of a $2n$-dimensional p-stratifold $\f S$
%$\{x_1,\dots,x_l\}$
having $(n-2)$-connected links of isolated singularities. Assume
that $n \equiv 6\; (\mod 8)$and that $\lbar{U} \hat \Sigma_i$ and
$\lbar{U} \hat \Sigma'_i$ are parallelizable with compatible
parallelizations on the boundary. Let further $e(\lbar{U} \hat
\Sigma_i)= e(\lbar{U} \hat \Sigma'_i)$ and $\sign(\lbar{U} \hat
\Sigma_i \cup_\Par \lbar{U} \hat \Sigma'_i) =0$. Then there is a
$k\in \{0,1\}$ such that $\hat{\f S} \sharp k(S^{n} \times S^{n})$
is almost equivalent to $\hat{\f S}'\sharp k(S^{n} \times S^{n})$.
\end{cor}

%--------------------------------------------------------------------------------
\section{4-dimensional results}
In this section we consider the exceptional case of a
$4$-dimensional p-stratifold and give a similar classification
result in that
situation. The proof of the main theorem can be found in \S \ref{proofs_4}.\\

For  a $4$-dimensional stratifold $\f S$,  every link of the
singularity $L_i$ is a 3-dimensional manifold. According to the
computation of $\Omega_*$ by Thom \cite{Th} we immediately obtain
from Theorem \ref{k}:
\begin{cor}
A four-dimensional p-stratifold with isolated singularities always
admits a resolution.
\end{cor}

If we further assume the links to be oriented we can use the
following well-known result, which can be proved easily.
\begin{prop}
Every orientable 3-manifold is parallelizable, hence in particular
spin.
\end{prop}

The normal 1-type of a simply connected 4-dimensional
spin-manifold is given by $B{\text{Spin}}$. Since
$\Omega^{\text{spin}}_3 = 0$ (see \cite[Lem. 9]{Mi_sp}) we obtain
the following corollary from Theorem \ref{g2}.
\begin{cor}
A four-dimensional p-stratifold with isolated singularities admits
an optimal resolution if and only if all links of singularities
are orientable.
\end{cor}

%Let us concentrate on resolutions by
%spin manifolds.\\

%\Def An (optimal) resolution $p:\hat{\f S} \lra \f S$ is called a
%{\em
%  spin-resolution} if $\lbar{U}\hat \Sigma_i$ is spin for a
%suitable representative of the neighbourhood germ.\\

%\Cl Since $\Omega^{\text{spin}}_3=0$ every 3-manifold bounds a spin
 %manifold.

We have to develop some notation in the topological category. Use
$B\text{TOP}$ to denote the classifying space of topological
vector bundles and let $B\text{TOPSpin}$ be the 2-connected cover
over $B\text{TOP}$. Let $M$ be a simply connected 4-manifold.
Using the Wu-Formula we can explain the Stiefel-Whitney-classes of
$M$. We call the topological manifold $M$ {\em spin} if $w_2(M)$
vanishes. One can show that the topological Gauss map of $M$ lifts
to $B\text{TOPSpin}$ if and only if $M$ is spin.  Note further
that if
such a lift exists, it is unique.\\

Using  \cite[Thm. 2]{Kr1} and the h-cobordism-Theorem in dimension
4 \cite[Thm. 10.3]{F} we formulate:

\begin{thm}\label{kr4}
Let $ M_1$  and $M_2$ be compact $4$-dimensional topological spin
manifolds with $e(M_1) =e(M_2)$ and  let $g: \Par M_1 \lra \Par
M_2$ be a homeomorphism compatible with the induced
spin-structures on the boundaries. If $M_1 \cup_g M_2$ vanishes in
$\Omega^{B\text{TOPSpin}}_4$, then $g$ can be extended to a
homeomorphism $G: M_1 \sharp k(S^2 \times S^2) \lra M_2\sharp
k(S^2 \times S^2)$ for $k \in \{0,1\}$.
\end{thm}

We call two resolutions {\em topologically equivalent} if the
diffeomorphism $\phi_{\b c}$ in the definition of equivalent
resolutions in \S \ref{int} is replaced by a homeomorphism. Using
this notation, we obtain the following classification result in
dimension four:

\begin{thm}\label{mt4}
Let $\hat{\f S} \lra \f S$ and
  $\hat{\f S}' \lra \f S$ be two
optimal resolutions of a $4$-dimensional p-stratifold $\f S$ with
isolated singularities $\{x_i\}_{i\in I}$, such that each link
$L_i$ is connected. Assume that both $\hat{\f S}$ and  $\hat{\f
S}' $ are spin and that for a suitable representative
$\lbar{U}\Sigma$ of the neighbourhood germ, the following
conditions hold for all $i \in I$:
\begin{itemize}
\item $e(\lbar{U}\hat \Sigma_i) = e(\lbar{U}\hat \Sigma'_i)$,
\item the spin-structures of $\lbar{U}\hat \Sigma_i$ and
$\lbar{U}\hat
  \Sigma'_i$ coincide on the boundary,
\item $\sign(\lbar{U}\hat \Sigma_i \cup_\Par \lbar{U}\hat
\Sigma'_i) =
  0$.
%\item $\lbar{U}\hat \Sigma_i \cup_{\Par}\lbar{U}\hat \Sigma'_i$ is elementary.
\end{itemize}
Then $\hat{\f S}\sharp k(S^2\times S^2)$ is topologically
equivalent to $\hat{\f S}'\sharp k(S^2\times S^2) $ for a $k \in
\{0,1\}$.
\end{thm}

%---------------------------------------------------------------------------
%\pagebreak
\section{Outline of the proofs}\label{proofs}

\subsection{Proof of Theorem \ref{k}}
Although the proof can be found in \cite{Kr3}, it is useful to
understand its nature for the succeeding results.

One of the basic tools for  constructing  a resolving map is the
following lemma, which can be proved with the help of Morse
theory, cf.\ \cite{Kr3}.
\begin{lem}\label{ml}
Let $W$ be a smooth compact manifold with boundary. Then there is a codense
compact subspace $X$ of $W$
 and a continuous map $f: \Par W \lra X$
such that $W$ is homeomorphic to $\Par W \times [0,1] \cup_f X$, where on $\Par W \times [0,1)$ the homeomorphism can be chosen to
be a diffeomorphism.
\end{lem}
In other words, every smooth manifold with boundary arises from
its collar by attaching a codense set. 
The notation {\em codense} stands for the complement of a dense subset. 
With this information we are
ready to
prove Theorem \ref{k}.\\

\begin{proof}
Let $p: \hat{\f S}  \lra \f S$ be a resolution. Set $W_i :=
\lbar{U}\hat\Sigma_i$. Since $W_i$ is a compact manifold with
boundary $L_i$,
we obtain $[L_i] = 0$ in $\Omega_{n-1}$.\\

Let on the other hand $W_i$ be a compact manifold bounding $L_i$
and let $f(\inter{N})$ be the top stratum of $\f S$. Set $\hat{\f
S} := N\cup (\sqcup_i W_i)$ and construct with the help of the
last
 lemma the following resolving map:
\begin{displaymath}
{\xymatrix{ {\hat{\f S} \cong  N  \cup_i(\Par W_i \times
[0,\varepsilon_i] \cup_{f_i} X_i)}  \ar[r]^p & { N \cup_i (L_i
\times [0,\varepsilon_i]  \cup\{x_i\}) \cong \f S} }}.
\end{displaymath}
\vspace{9mm}
\hspace{0.1cm}
\begin{picture}(12,-1)
\put(4.3,0.4){\oval(5.6,1)[b]}
\put(1.5,0.4){\line(0,1){0.2}}
\put(7.1,0.4){\vector(0,1){0.2}}
\put(4.6,0){\small id}

\put(6.0,0){\oval(5.2,1)[b]}
\put(3.4,0){\line(0,1){0.5}}
\put(8.6,0){\vector(0,1){0.5}}
\put(6.3,-0.45){\small id}

\put(7.95,-0.3){\oval(5.1,1)[b]}
\put(5.4,-0.3){\line(0,1){0.8}}
\put(10.5,-0.3){\vector(0,1){0.8}}
\end{picture}

\vspace{-10mm}

\end{proof}

\subsection{Proof of Theorem \ref{nt}}
We consider a $2n$-dimensional manifold $M$, which is
$(n-1)$-connected, and want to determine its $(n-1)$ type. We
begin with the classification up to
fiber homotopy equivalence of fibrations $p: B \lra BO$, with a CW-complex $B$,
fulfilling
\begin{enumerate}
\item $B$ is $(n-1)$-connected and
\item $\pi_i(F) = 0$ for $i \geq n$, where
$F$ is the fiber.
\end{enumerate}
Compare such a fibration with
the $(n-1)$-connected cover of $BO$:
\begin{displaymath}
\xymatrix {
         & BO\langle n-1 \rangle  \ar[d]^{p_{\langle n-1\rangle }}
 & F_{\langle n-1\rangle}  \ar@{_{(}->}[l] \\
B \ar[r]_p \ar@{-->}[ur] & BO
}
\end{displaymath}
Since all obstructions vanish, we obtain a lift
$\bar p: B \lra BO\langle n-1\rangle$, which without loss of generality may be assumed
to be a fibration. From the long exact homotopy sequence we see that
the homotopy groups of the fiber vanish, except in dimension $(n-1)$,
where the group is $\pi := \coker (p_* :\pi_n(B) \lra \pi_n(BO))$.
Thus, $\bar p: B \lra BO\langle n-1\rangle$ is a fibration with fiber
$ K(\pi,n-1)$. Such fibrations are classified in \cite[\S 5.2]{Ba} as follows:
\begin{displaymath}
\begin{array}{ccc}
      [BO\langle n-1\rangle,K(\pi,n)]/[Aut(\pi)] & \stackrel{\cong}{\lra} &  \mathcal{F}(K(\pi,n-1),BO\langle n-1\rangle) \\
             {[f]} & \longmapsto & f^*(P(K(\pi,n))
\end{array}
\end{displaymath}
Here
$\mathcal{F}(K(\pi,n-1),BO\langle n-1\rangle)$ denotes the set of all fibrations over $BO\langle n-1\rangle$
with fiber $K(\pi,n-1)$ up to fiber homotopy equivalence.
Thus $\bar p: B \lra BO\langle n-1\rangle$ is a pull back
$$
\xymatrix{
B \ar[r] \ar[d]^{\bar{p}} & P(K(\pi,n)) \ar[d] \\
BO\langle n-1\rangle  \ar[r]_{\theta} & K(\pi,n)
}
$$
with an appropriate map $\theta$. The definition of $\pi$
force the induced map
$\theta_*:\pi_n(BO\langle n-1\rangle ) \lra \pi_n(K(\pi,n))=\pi$
to be surjective. Therefore we can assume $\theta_*$ to be the canonical
projection to the factor group $\pi$. On the other hand, each factor
group of $\pi_n(BO)$ leads to a fibration with the claimed properties.
We summarize this discussion in
\begin{lem} \label{nt1}
Fibrations with properties 1. and 2. are given by
\mbox{$p_k: B_k \lra BO$} up to fiber homotopy equivalence
($k \in \Bbb N, 0 \leq\ k <|\pi_n(BO)|$).
\end{lem}

Consider now the fibration  $p: B_k \lra BO$ and ask for a lift:
\begin{displaymath}
\xymatrix{& B_{k} \ar[d]^{p} \\
 \nu: M \ar[r] \ar@{-->}[ur] & BO }
\end{displaymath}
With the help of obstruction theory we see that such a lift exists
if and only if

\begin{displaymath}
  \im(\nu_*:\pi_n(M) \lra \pi_n(BO)) \leq \langle kx \rangle,
\end{displaymath}
 where $\pi_n(BO) = \langle x \rangle$.
Combining this discussion with Lemma \ref{nt1}, the statement of
Theorem \ref{nt}
follows immediately.

%---------------------------------------------------------------------
\subsection{Surgery in the middle dimension}
First we give a brief introduction in surgery, for more details
compare \cite{W1}.

Surgery is a tool to eliminate homotopy classes in the category of
manifolds. Let $M$ be a compact $m$-dimensional manifold. One starts
with an embedding
$f: S^r \times D^{m-r} \hookrightarrow \inter{M}$ and define
$T := D^{r+1} \times D^{m-r} \cup_f (M \times I)$,
where $f$ is considered as a map to $M \times {1}$. The corners of the
manifold $T$ can always be straighten, according to \cite{C-F}. This
construction is called {\em attaching an $(r+1)$-handle} and $T$ the
{\em trace of a surgery via $f$}.

The boundary of $T$ is $M \cup (\Par M \times I) \cup M'$ and we call
$M'$ the {\em result of a surgery of index $r+1$ via $f$}. It is not
difficult to see that $T$ can also be viewed as the trace of a surgery
on $M'$ via the obvious embedding of $D^{r+1}\times S^{m-r-1}$ into
$M'$, compare \cite{Mi0}.\\

Since we are working in the category of manifolds with $B$-structures
we have to ask, whether the result of surgery via an embedding $f$ is
equipped with a $B$-structure. For general results see \cite{Kr1}. In
our situation we are only looking at fibrations $p_k:B_k \lra B$ defined
in \S \ref{cl}.
\begin{lem}\label{bs}
Let $M$ be a manifold of dimension $m \in \{2n, 2n+1\}$ with  $B_k$-structure and $f: S^n \times
D^{m-n} \hookrightarrow M$ an embedding. Then $\bar \nu: M \lra
B$ extends to a normal $B_k$-structure on $T$, the trace of the
surgery via $f$.
\end{lem}
\begin{proof}
The embedding $f: S^n \times D^{m-n} \hookrightarrow M$ induces a normal
$B_k$-structure on $S^n \times D^{m-n}$ denoted by $f^*\bar
\nu$. There is a unique (up to homotopy) $B_k$-structure on $D^{n+1}
\times D^{m-n}$ and we have to show that its restriction to $S^n
\times D^{m-n}$ is $f^* \bar \nu$. Let $F_k$ be the fiber of $p_k : B_k
\lra BO$. From the long exact homotopy sequence
we see that the different $B_k$ structures on $S^n \times D^{m-n}$ are
classified by $\pi_n(F_k)$. But the construction of $B_k$ implies that
$\pi_n(F_k)=0$, thus the restriction on $S^n \times D^{m-n}$ coincides with the given
structure and we obtain a  $B_k$-structure on $T$.
\end{proof}
\begin{lem}\label{sur1}
Let $M$ be an $(n-1)$-connected $2n$-dimensional manifold. Let
$H_n(M)$ have a free basis $\{\lambda_1, \dots,
\lambda_r,\mu_1,\dots,\mu_r\}$ with
\begin{itemize}
\item $\lambda_i \cdot \lambda_j = 0$ for all $i,j$ and
\item $\lambda_i \cdot \mu_j = \delta_{ij}$ for all $i,j$.
\end{itemize}
If further each embedded sphere representing a generator $\lambda_i$
has a trivial normal bundle, then by a finite
sequence of surgeries the homology group $H_n(M)$ can be eliminated.
\end{lem}
\Rem The last two lemmas imply that the condition in Theorem \ref{hs} is sufficient.

\begin{proof}
According to a result of Haefliger ({\cite{Hae}}), every element of
$H_n(M)$ can be represented by an embedding $S^n \hookrightarrow M$.
Let $\lambda :=\lambda_r$ and $\varphi_0: S^n\hookrightarrow M$ an
embedding representing $\lambda$. Since $\lambda$ is assumed to have a
trivial normal bundle, the embedding can be extended to
$\varphi: S^n \times D^n \hookrightarrow M$. Set  $M_0 := M
-\varphi((S^n \times D^n)^{\circ})$ and   let $M'$ denote the result
of surgery via $\varphi$, i.e.\
$M' = M_0 \cup_{\varphi} (D^{n+1} \times S^{n-1})$.
%For explicit explanations on surgery compare \cite{W1}.
Combine now
the long exact homology sequences of pairs $(M,M_0)$ and $(M',M_0)$
and obtain the following commutative diagram:
\begin{displaymath}
  \xymatrix{
     & \mathbb Z \ar[d]^{\tilde{\Par '}} \ar[dr]^{1 \mapsto \lambda} \\
0\ar[r] & H_n(M_0) \ar[r]_{i_*}\ar[d] & H_n(M) \ar[r]^(.6){\cdot \lambda} &
  \mathbb Z \ar[r]^<(.2){\tilde{\Par}} & H_{n-1}(M_0) \ar[r] & 0 \\
     & H_n(M') \ar[d] \\
     & 0 \\
   }
\end{displaymath}
The excision, together with the surjectivity of $H_n(M)\stackrel{
  \lambda \cdot}{\lra} \mathbb Z$, implies $M_0$ is $(n-1)$-connected,
further we see
$ H_n(M_0) =
\langle \lambda_1,\dots,\lambda_r,\beta_1,\dots,\beta_{r-1}\rangle$.

Thus $M'$ is $(n-1)$-connected as well and
\begin{displaymath}
  H_n(M') = \langle  \bar{\lambda}_1,\dots,\bar{\lambda}_{r-1},
                         \bar{\mu}_1,\dots,\bar{\mu}_{r-1} \rangle,
\end{displaymath}
where the generators are given by $\bar{\lambda}_i=\lambda_i + \lambda \mathbb Z$ and
$ \bar{\mu}_i = \mu_i + \lambda \mathbb Z$. We can always deform the embedding
of the generator $\bar \lambda_i$ to $M_0$, such that it represents
the class ${\lambda}_i + {\lambda}_r \in H_n(M_0)$. Thus we conclude
$
%\begin{displaymath}
\begin{array}{lcl}
  \bar{\lambda}_i\cdot\bar{\lambda}_j = 0 & \mbox{und} &
    \bar{\lambda}_i\cdot\bar{\mu }_j = \delta_{ij}
\end{array}
$.
%\end{displaymath}
 Since the
intersection product $\lambda_i \cdot \lambda_r$ vanishes we obtain
$\nu({\lambda}_i + {\lambda}_r) = \nu({\lambda}_i)+\nu({\lambda}_r)$
(cf.\ \cite{W0}). Now we proceed with the manifold $M'$ and
inductively obtain the desired statement.
\end{proof}

\subsection{Surgery on odd-dimensional manifolds}
\begin{lem}{\label{sur2}}
Let $T$ be a bordism in $\Omega_{2n}(B_k)$ between a manifold
$M$ of index $k$ and a
homotopy sphere $\mathcal S$. Then $T$ is bordant in $\Omega_{2n}(B_k)$
rel. boundary to $T'$, such that its homology groups
 $H_n(T')$ and $H_{n+1}(T')$
are free and $H_i(T') = 0$ for $i \not\in \{0,n,n+1,2n+1\}$.
\end{lem}
\begin{proof}
According to Theorem \ref{k2}, we can assume that
$T$ is $(n-1)$-connected, further the Universal Coefficient Theorem
implies that $H_{n+1} (T)$ is free and $\mbox{Tor}(H_n(T)) \cong
\mbox{Tor}(H_n(T,M))$. We will show that the torsion of $H_n(T,M)$ can
be eliminated by a finite sequence of surgeries.
%Hence we turn our attention on
%the $n$-th homology group.

From the long exact homology sequence of the pair $(T,M)$ we see, that
every torsion element $\alpha' \in H_n(T,M)$ comes from an element
$\alpha \in H_n(T)$. After possible correction of  $\alpha$ by an element of $H_n(M)
\cong \pi_n(B_k)$ we achieve $\nu (\alpha) = 0$.\\
  Let
$\varphi: S^n \times D^{n+1} \hookrightarrow \inter T$ be an embedding
representing $\alpha$. As in the previous proof, we set
 $T_0 = T - \varphi(( S^n \times D^{n+1})^{\circ})$ and $T'= T_0 \cup
 (D^{n+1} \times S^n)$. We combine now the exact triple sequences for
$(T,T_0,M)$ and $(T',T_0,M)$ to obtain:\newpage

\noindent\hbox{}\vglue -10mm
\begin{displaymath}
  \xymatrix{
& & 0 \ar[d] \\
& & H_{n+1}(T_0,M)\ar[d] \\
& & H_{n+1}(T',M) \ar[d]^{\beta \cdot}\\
& & \Bbb Z \ar@{|->}[dr]^{1\mapsto \alpha '}\ar[d]^{d'} \\
H_{n+1}(T,M)\ar[r]^(0.7){\alpha \cdot} &
 \Bbb Z \ar[r]^<(.2)d\ar@{|->}[dr]^{1\mapsto \beta '} & H_n(T_0,M)\ar[r] \ar[d]^{j'_*}
& H_n(T,M) \ar[r] & 0 \\
& & H_n(T',M) \ar[d]\\
& & 0 \\
}
\end{displaymath}
The element $\beta \in H_n(T')$ is given by the embedding $\psi:
D^{n+1} \times S^n \hookrightarrow T'$ and corresponds to the
homotopy class $[\psi|_{0 \times S^n}]$.

We consider the two cases, where $\alpha$ is free or a torsion element
modulo \linebreak$i_*(H_n(\Par T))$, separately.

\medskip
{\bf Case 1}\qua $\alpha$ is primitive $(\mod \; i_*(H_n(\Par T)))$.

In this case the Poincar\'{e} duality implies that the map
$H_{n+1}(T,M) \stackrel{\alpha \cdot}{\lra} \Bbb Z$ is surjective,
therefore \begin{displaymath}
  H_n(T',M) \cong H_n(T,M)/{<\alpha'>}.
\end{displaymath}
Hence the torsion group of $H_n(T',M)$ has been reduced.

\medskip
{\bf Case 2}\qua $\alpha$ is torsion $(\mod \; i_*(H_n(\Par T)))$.

The map $H_{n+1}(T,M) \stackrel{\alpha \cdot}{\lra} \Bbb Z$ is
trivial now.  %and we see from the sequence
Denote with $o(x)$ the order of a torsion element $x$. From the
sequence above we see that $o(\alpha')d(1) \subset \im d$, thus
there exists a $b' \in \Bbb Z$ such
that
\begin{equation*}
o(\alpha') d'(1) = b' d(1). \tag{$*$}
\end{equation*}

If $b'=0$, then the element $\beta'$, corresponding to $\alpha'$, has
infinite order, and the torsion rang again decreases.

If $b' \neq 0$, then $(\ker j'_*) \subset \mbox{Tor}(H_n(T_0,M))$, thus
$j'_*$ is injective on the free part of $H_n(T_0,M)$, therefore  the
element $d'(1)$
has infinite order and $o(\beta')\mid |b'|$. We need a finer case differentiation.\\

\medskip
{\bf Claim}\qua If $n$ is even and $\alpha$ a torsion element of
order $a$ in $H_n(T)$, then $\beta \in H_n(T')$ is an element of
infinite order. \\
Consider the pair sequences $(T,T_0)$ and $(T', T_0)$ and obtain the
following diagram
\begin{displaymath}
  \xymatrix{
 & & & & 0 \ar[d] \\
 & & & & H_{n+1}(T_0)\ar[d] \\
 & & & & H_{n+1}(T') \ar[d]\\
 & & & & \Bbb Z \ar@{|->}[dr]^{1\mapsto \alpha }\ar[d]^{d'} \\
0 \ar[r] & H_{n+1}(T_0) \ar[r] & H_{n+1}(T)\ar[r]^{\alpha \cdot} &
 \Bbb Z \ar[r]^d\ar@{|->}[dr]^{1\mapsto \beta } & H_n(T_0)\ar[r] \ar[d]
& H_n(T) \ar[r] & 0 \\
 & & & & H_n(T') \ar[d]\\
 & & & & 0 \\
}
\end{displaymath}
As in the previous case,
there exists  a $b \in \Bbb Z$ such that
\begin{displaymath}
        a\cdot d'(1) = b\cdot d(1) \;\; \mbox{in} \;H_n(T_0)
\end{displaymath}
or equivalently after the identification of the generators
\begin{displaymath}
  (\varphi_0)_*(a\cdot([S^n]\otimes 1) - b\cdot(1\otimes[S^n] )) = 0.
\end{displaymath}
By computing the self intersection number of $a\cdot([S^n]\otimes
1) - b\cdot(1\otimes[S^n] ) $ we obtain $2ab([S^n]* \otimes
[S^n]^*) \in \im ((\varphi_0)^*: H^{2n}(W_0) \lra H^{2n}(S^n\times
S^n))$. Thus $2ab = 0$, and since $a \neq 0$ it follows $b = 0$.
We conclude again that $\beta$ is of infinite order, as described
previously.\\

\medskip
{\bf Claim}\qua If $n$ is even and $\alpha$ has infinite order,
then the element $\alpha'$ will be eliminated.\\
In this situation we have $\alpha = \tilde{\alpha} + i_*({\gamma})$, where $\tilde{\alpha} $
is a torsion element and ${\gamma} \in H_n(M)$. The normal bundle
map is given by
\begin{displaymath}
  \nu: H_n(T) \lra \pi_{n-1}(SO({n+1})) = \pi_{n-1}(SO).
\end{displaymath}
If $\pi_{n-1}(SO) = 0$ or $\Bbb Z$ the fact that the map $\nu: H_n(T) \lra
\pi_{n-1}(SO)$ is a homomorphism implies, that $\tilde{\alpha}$
already has a trivial normal bundle, this leads us in the situation of
the last claim.
It remains to study the situation $\pi_{n-1}(SO) = \Bbb Z/2$ with
$\nu(\tilde\alpha) = 1$. Observe that without loss of generality we can assume
$i_*(\gamma)$ to be primitive.
 Thus the map $H_{n+1}(T) \stackrel{\alpha
  \cdot}{\lra} \Bbb Z$ from the sequence of the last claim is
surjective and we obtain
\begin{displaymath}
  H_n(T') \cong H_n(T)/<\alpha>  \;\;\Rightarrow\;\;
  H_n(T',M) \cong H_n(T,M)/<\alpha'>.
\end{displaymath}

{\bf Claim}\qua If $n$ is odd then the torsion group
$\mbox{Tor}(H_n(T ,M))$ can be reduced.\\
The group $\pi_{n}(SO_{n+1})$ acts on the trivializations via
\begin{displaymath}
\begin{array}{lcllcl}
  (\omega,\varphi) &\longmapsto& \psi : & S^n \times D^{n+1} &\hookrightarrow & T \\
 & & & (x,y) & \mapsto &  \varphi(x,\omega(x)y) \\
\end{array}
\end{displaymath}
Note, that the change of trivialization does not affect the
$B$-structure on $S^n \times D^{n+1}$ since the induced $B$-structures
$\varphi^*\bar \nu$ and $\psi^*\bar\nu$ differ by an element from
$\pi_n(F_k) = 0$.\\
Let $y_0$ be the basis point of $S^n$, then we see
$$
\psi_0(x,y_0) = \varphi_0(x,\omega(x)(y_0)) = \varphi_0(x,pw(x)) =
\varphi_0(\id \times p\omega)\Delta(x),
$$
where the map $p : SO(n+1) \lra S^n$ is the canonical fiber bundle and
$\Delta : S^n \lra S^n \times S^n$ the diagonal. Denote by
$i_l : S^n \lra S^n \times S^n$ the inclusion in the $l$-th component
and
let $\iota$ be the generator of $\pi_n(S^n)$, set $\iota_l:=(i_l)_*(\iota)$. We pass to homotopy and
use the fact that the map $\pi_* \delta : \pi_{n+1}(S^{n+1})\lra
\pi_n(S^n)$ from the long exact homotopy sequence for $p$ is given by
multiplication with $2$ \cite{Ste}. Thus we compute
\begin{displaymath}
\begin{array}{lcl}
  (\psi_0)_*(\iota_1) & = & (\varphi_0)_*(\id \times p\omega)_*\Delta_*(\iota) \\
                             & = & (\varphi_0)_*(\id \times p\omega)_*(\iota_1 +\iota_2 ) \\
                             & = & (\varphi_0)_*( \iota_1 + 2k\cdot \iota_2) \\
%& = & (?????????????????????????????????????) \\
                             & = & (\varphi_0)_*(\iota_1) + 2k\cdot (\varphi_0)_*(\iota_2)
\end{array}
\end{displaymath}
Since $M$ is $(n-1)$-connected we obtain the corresponding statement in
homology:
\begin{displaymath}
(\psi_0)_*([S^n]\otimes 1) = (\varphi_0)_*([S^n]\otimes 1) + 2k\cdot (\varphi_0)_*(1\otimes [S^n])
\end{displaymath}
%\end{displaymath}
The equality $(*)$ now becomes
\begin{displaymath}
\begin{array}{l}
  b \cdot (\psi_0)_*(1\otimes [S^n]) = a \cdot ((\psi_0)_*([S^n]\otimes 1) -
2k(\psi_0)_*(1\otimes [S^n])) \\
a \cdot (\psi_0)_*([S^n]\otimes 1) = (b + 2ka) (\psi_0)_*(1 \otimes [S^n])
\end{array}
\end{displaymath}
Since the case $b+2ka=0$ has already been treated, we assume   $b +2ka
\neq 0$. By choosing $k$ appropriately we can achieve
\begin{displaymath}
  o(\beta ') \leq o(\alpha ').
\end{displaymath}
Let $p$ be a prime number such that $(\alpha ')_p \neq 0$ in
$H_n(T,M;\Bbb F _p)$.
From the analogous sequences with $\Bbb F_p$-coefficients we conclude
\begin{displaymath}
H_n(T ',M;\Bbb F _p) \cong H_n(T_0,M;\Bbb F _p)/{\mbox{im} \, d'} \cong
H_n(T,M;\Bbb F _p)/{<\alpha '>}
\end{displaymath}
and with universal coefficient theorem
\begin{displaymath}
  |\mbox{Tor}(H_n(T ',M))| < |\mbox{Tor}(H_n(T ,M))|.
\end{displaymath}
Combining all cases together we see, that a torsion element of
$H_n(T,M)$ can either be eliminated by a finite sequence of surgeries
or can be replaced by an element of infinite order. Inductively we
obtain the desired statement.
\end{proof}
%------------------------------------------------------------------------------------
\subsection{Proof of Theorem \ref{hs}}
\begin{proof}
We only have to prove that every manifold normally $B_k$ bordant to a
homology sphere is elementary. \\ Let $T$ be a bordism between $M$ and
a homology sphere $\mathcal S$. According to Theorem \ref{k2}, we can
without loss of generality assume that $T$ is $(n-1)$-connected and
that the homology groups $H_n(T)$ and $H_{n+1}(T)$ are free. The long
exact sequence of the pair $(T,M)$ together with Poincar\'{e} duality
leads to the following commutative diagram: {
\fontsize{9}{15} \selectfont
\begin{displaymath}
  \xymatrix{
%0 \ar[r] &
H_{n+1}(T)\ar[r] \ar[d]^{P.-D.} & H_{n+1}(T,M) \ar[r]^{\delta}\ar[d]^{P.-D.}
 &
H_n(M) \ar[r]^{i_*} \ar[d]^{P.-D.}&
 H_n(T) \ar[r]^{j_*} \ar[d]^{P.-D.} & H_n(T,M) \ar[r]  & 0 \\
% &
H^n(T,M) \ar[r]^{j^*}  &  H^n(T) \ar[r]^{i^*}&  H^n(M) \ar[r]^d &  H^{n+1}(T,M) \\
}
\end{displaymath}
}
From this we see \begin{displaymath}
  \dim H_n(M) = 2(\dim  \ker i_*).
\end{displaymath}
Thus $\ker i_*$ is a direct summand of $H_n(M)$. It is not hard to
verify, that  $\ker i_*$ is isotropic. Let now $\lambda: S^n \hookrightarrow M$
be a representative of an element of  $\ker i_*$, then the map can be
extended to $\bar\lambda: D^n \lra T$ and with the help of the Whitney-trick
\cite{Hae} this map can be assumed to be an embedding. Thus the
restriction of the normal bundle of $D^n$ to the boundary gives us a
trivialization of the normal bundle of $\lambda$.
\end{proof}

%---------------------------------------------------------------------------------
\subsection{Proof of Theorem \ref{mt}}\label{pmt}
\begin{proof}
Using the statement of Theorem \ref{hs} we see, that the
conditions of Theorem \ref{k2} are fulfilled for the manifolds
$\lbar{U}\hat \Sigma_i \sharp \mathcal S$ and $\lbar{U}\hat
\Sigma'_i$ for a homotopy sphere $\mathcal S$. Thus, the
diffeomorphism on the boundary can be extended to the the
diffeomorphism $\lbar{U}\hat \Sigma_i \sharp \mathcal S \sharp
k(S^n\times S^n) \lra \lbar{U}\hat \Sigma'_i\sharp k(S^n\times
S^n) $ for $k\in\{0,1\}$. From the definition of resolutions and
from construction the diffeomorphism on the boundaries we see that
the obtained diffeomorphisms extend in the obvious way to a
diffeomorphism $\f S \sharp \mathcal S \sharp k(S^n\times S^n)
\lra \f S\sharp k(S^n\times S^n)$ having
the desired properties. \\

It remains to show, that the conditions of the theorem are true
for every pair of representatives of the neighbourhood germs
$[\lbar{U}\hat \Sigma_i]$ and $[\lbar{U}\hat\Sigma'_i]$, once we
have checked them on  a single representative pair. If ${\b c}_i:
L_i \times [0,\varepsilon_i/2]\lra N $ and ${\b d}_i: L_i \times
[0,\varepsilon'_i/2]\lra N $ with $\varepsilon_i < \varepsilon'_i$
are two representatives of the germ of collars around $L_i$, then
there exists a $\delta_i > 0$ such that  ${\b c}_i$ coincides with
${\b d}'_i$ on $ L_i \times [0,\delta_i]$.  We choose a
diffeomorphism $\eta_i:[0,\varepsilon_i/2] \lra
[0,\varepsilon'_i/2]$ with $\eta_i\equiv \id$ on $[0,\delta_i]$.
The map induces an isomorphism
$$
\begin{array}{lcl}
\lbar{U} {\Sigma_i}_{{\b c}} &\lra& \lbar{U} {\Sigma_i}_{{\b d}}\\
f({\b c}_i(x,t)) &\longmapsto& f({\b d}_i(x,\eta(t)))
\end{array}
$$
being the identity on a small neighbourhood of $x_i \in \Sigma$.
This gives us a diffeomorphism between
${\lbar{U}{{\hat{\Sigma}_i}}}_{{\b c}}$ and
${\lbar{U}{\hat\Sigma_i}}_{{\b d}}$ making the following diagram
commutative:
$$
\xymatrix{ \Par {\lbar{U}{{\hat{\Sigma}_i}}}_{{\b c}}
\ar@{^(->}[r] \ar[d]^{\cong} & {\lbar{U}{{\hat{\Sigma}_i}}}_{{\b
c}} \ar[r]^{\cong} & {\lbar{U}{{\hat{\Sigma}_i}}}_{{\b
d}}\ar@{<-^)}[r] &
\Par {\lbar{U}{{\hat{\Sigma}_i}}}_{{\b d}}\ar[d]^{\cong} \\
%--------------------------------
\Par {\lbar{U}{{\hat{\Sigma'}_i}}}_{{\b c}} \ar@{^(->}[r] &
{\lbar{U}{{\hat{\Sigma'}_i}}}_{{\b c}} \ar[r]^{\cong} &
{\lbar{U}{{\hat{\Sigma'}_i}}}_{{\b d}}\ar@{<-^)}[r] &
\Par{\lbar{U}{{\hat{\Sigma'}_i}}}_{{\b d}} }
$$
This completes the proof.
\end{proof}

%-------------------------------------------------------------------------------
\subsection{Algebraic invariants}\label{proofs_alg}

ad (2)\qua Let $(H, \Lambda, s\nu_*)$ be elementary and ${\cal L}=
\langle\lambda_1,\dots \lambda_k\rangle$ a Lagrangian with
$s\nu_*|_{\cal L} \equiv 0$. Thus $\sign(\Lambda)=0$ and $0 =
s\nu_*(\lambda_i) =
\Lambda\langle\kappa_{s\nu_*},\lambda_i\rangle$ for all $1\leq
i\leq k$. Since ${\cal L}$ is maximal it follows that
$\kappa_{s\nu_*} \in {\cal L}$ and therefore
$s\nu_*(\kappa_{s\nu_*}) =0$.

On the other hand let $\sign(\Lambda) = 0$ and
$s\nu_*(\kappa_{s\nu_*}) =0$. Choose a basis $\{
\lambda_1,\dots,\lambda_k,\mu_1,\dots,\mu_k\}$ of $H$ such that $
\Lambda(\lambda_i,\lambda_j)=0$  and  $\Lambda(\lambda_i,\mu_j) =
\delta_{ij}$.

There is nothing to show if $\kappa_{s\nu_*} = 0$. Otherwise
 we can without loss of generality
assume that $s\nu_* (\lambda_i) = 0$ for all $i > 1$, since
$s\nu_*$ is a homomorphism.
%and $s\nu_* (\mu_j) = 0$ for all $j > 2$.
%For $\lambda : H \lra H^*$ is an
%isomorphism there is an $x = \sum_{i=0}^k(a_i \lambda_i+b_i\mu_i) \in
%H$, %$k_i,l_i \in \B Z$,
%such that
Recall the equality $
  s\nu_*(v) = \Lambda(\kappa_{s\nu_*},v) \;\;\forall v \in H
$. Since $\kappa_{s\nu_*} \in H$ there are $a_i,b_i \in \Bbb Z$
such that $\kappa_{s\nu_*} = \sum_{i=1}^k(a_i
\lambda_i+b_i\mu_i)$. Consider a sub-Lagrangian ${\cal L}' :=
\langle \lambda_2,\dots,\lambda_k\rangle$. If $\kappa_{s\nu_*}
\not\in {\cal L}'$, build ${\cal L} :=  \langle
\lambda_2,\dots,\lambda_k,\tilde{\kappa}_{s\nu_*}\rangle$, where
$\tilde{\kappa}_{s\nu_*}$ is a primitive element of $H$ with
$\kappa_{s\nu_*} \in \langle\tilde{\kappa}_{s\nu_*}\rangle$. This
is a Lagrangian, satisfying $s\nu_*|_{\cal L} \equiv 0$. In the
case of $\kappa_{s\nu_*} \in {\cal L}'$, the coefficients
$a_1,b_1,\dots,b_k$ have to be zero, thus ${\cal L} = \langle
\lambda_1,\dots,\lambda_k\rangle$ is a Lagrangian with the
desired property. \\

ad (3)\qua The conditions are obviously necessary. To see that they
are also sufficient choose a symplectic basis $\{\lambda_1,
\dots,\lambda_k,\mu_1,\dots,\mu_k\}$ of $H$. Sort the generators
in the following way
\begin{displaymath}
  \begin{array}{lcr}
s\nu_*(\lambda_i) = s\nu_*(\mu_i)=1 & \hspace{2cm} & \mbox{for } i\leq s, \\
s\nu_*(\lambda_i) = 0 & & \mbox{for } i > s,
  \end{array}
\end{displaymath}
where $s$ is an integer between $0$ and $k$. The assumption 
$$\Phi(H) = \sum_{i=1}^k s\nu_*(\lambda_i)s\nu_*(\mu_i) = 0$$ implies
that $s \equiv 0 (\mod 2)$. Construct a new basis
$\{\lambda'_{1},\dots,\mu'_k\}$ for $H$ by the substitution
\begin{displaymath}
  \begin{array}{lcl}
\lambda'_{2i-1} = \lambda_{2i-1} + \lambda_{2i}, & \mbox{  } &
\lambda'_{2i} = \mu_{2i-1} - \mu_{2i}, \\
\mu'_{2i-1}=\mu_{2i} & \mbox{  } & \mu'_{2i} = \lambda_{2i}
  \end{array}
\end{displaymath}
for $2i \leq s$, and
\begin{displaymath}
\begin{array}{lcl}
\lambda'_i = \lambda_i & \mbox{  } & \mu'_i = \mu_i
\end{array}
\end{displaymath}
for $i>s$. This new basis is again symplectic and satisfies the
condition
\begin{displaymath}
  s\nu_*(\lambda'_1) = \dots =  s\nu_*(\lambda'_k) = 0.
\end{displaymath}

%--------------------------------------------------------------------------
\subsection{Proof of Theorem \ref{mt4}}\label{proofs_4}

\begin{proof}
As in the proof of Theorem \ref{mt} we conclude that it is enough
to show that the diffeomorphism on the boundary $\Par \lbar{U}\hat
\Sigma_i \lra  \Par \lbar{U}\hat \Sigma'_i$ can be extended to a
homeomorphism on $\lbar{U}\hat \Sigma_i \sharp k(S^2\times S^2)
\lra \lbar{U}\hat \Sigma'_i\sharp k(S^2\times S^2)$. Since the
resolutions are optimal, the manifolds $\lbar{U}\hat \Sigma_i$ and
$\lbar{U}\hat \Sigma'_i$ are $1$-connected, hence $M :=
\lbar{U}\hat \Sigma_i \cup_{\Par} \lbar{U}\hat \Sigma'_i$ is again
1-connected. In order to apply Theorem \ref{kr4}  we have to show
that the closed $4$-dimensional spin manifold with vanishing
signature is bordant to a homotopy sphere. Then we apply the
topological 4-dimensional Poincar\'e conjecture proved by Freedman
\cite[Thm. 1.6]{F} and obtain the desired statement.

We  want to use surgery to prove that $M$ is bordant to a homotopy
sphere. Since $M$ is spin and $\sign(M) =0$,  there is a basis
$\{\lambda_1,\dots,\lambda_k,\mu_1,\dots,\mu_k\}$ of $H_2(M)$
satisfying
$$
\Lambda(\lambda_i,\lambda_j) = 0 \;\; \Lambda(\mu_i,\mu_j)=0 \;\;
\Lambda(\lambda_i,\mu_j) = \delta_{ij}.
$$
We can not use the Haefliger's embedding theorem   in dimension 4,
but according to \cite[Thm. 3.1, 1.1]{F} every generator
$\lambda_i$ is represented by a topological embedding $S^2 \times
D^2 \hookrightarrow M$. Knowing this we can proceed in exactly the
same way as in the proof of Lemma \ref{sur1}, working in the
category TOP.
%Hence we can proceed as in the proof of Lemma \ref{sur1}, working
%in the category TOP.
%
%Lemma \ref{sur2} is still valid in dimension 4 as well as the
%arguments in \S \ref{pmt}.  \\
%
%To complete the proof, observe that according to arguments from \S
%\ref{pmt} it is enough to check
%the conditions for a single representative of the neighbourhood germ.

%We can not use the embedding Theorem of Haefliger in dimension 4,
%but according to \cite[Thm. 3.1, 1.1]{F} every generator
%$\lambda_i$ is represented by a topological embedding $S^2 \times
%D^2 \hookrightarrow M$. proof of Lemma \ref{sur1}, working in the
%category TOP.

Lemma \ref{sur2} is still valid in dimension 4 as well as the
arguments in \S \ref{pmt}.  \\

To complete the proof observe that according to arguments from \S
\ref{pmt} it is enough to check the conditions for a single
representative of the neighbourhood germ.
\end{proof}

\Addresses\recd 
\end{document}

%% file: agtout.tex
\def\ifplaintex{\expandafter\ifx\csname documentclass\endcsname\relax}

\def\gtp{{\mathsurround=0pt\it $\cal G\mskip-2mu$eometry \&\ 
$\cal T\!\!$opology $\cal P\!$ublications}}  % GT publications

\def\recd{{\small Received:\qua\receiveddate\ifx\reviseddate\relax
\else\qquad Revised:\qua\reviseddate\fi\par}} 

\def\lognumber#1{\def\thelognumber{#1}}
\def\volumenumber#1{\def\thevolumenumber{#1}}
\def\volumeyear#1{\def\thevolumeyear{#1}}
\def\papernumber#1{\def\thepapernumber{#1}}
\def\pagenumbers#1#2{\def\startpage{#1}\def\finishpage{#2}}
\def\published#1{\def\publishdate{#1}}

\def\received#1{\def\receiveddate{#1}}
\def\revised#1{\def\reviseddate{#1}}
\def\accepted#1{\def\accepteddate{#1}}

\long\def\asciiabstract#1{\long\def\theasciiabstract{#1}}

\let\\\par\let\thelognumber\relax\let\thevolumenumber\relax
\let\thepapernumber\relax\let\thevolumeyear\relax\let\startpage\relax
\let\finishpage\relax\let\publishdate\relax\let\receiveddate\relax
\let\reviseddate\relax\let\accepteddate\relax\let\theasciititle\relax
\let\theasciiauthors\relax
\let\theasciiabstract\relax

\let\theasciiemail\relax

\ifplaintex
\font\logobig=cmssbx10 scaled 3836
\font\logomed=cmssbx10 scaled 2557
\else
\font\logobig=cmssbx10 scaled 4200
\font\logomed=cmssbx10 scaled 2800
\fi

\long\def\makeagttitle{   %%% start of definition of \makeagttitle
\count0=\startpage
\agt\hfill      %   Journal title (top left) 
\hbox to 45truept{\vbox to 0pt{\vglue -13truept{\logomed A\kern -.37em{\logobig 
T}\kern -.38em G}\vss}\hss}
\break
{\small Volume \thevolumenumber\ (\thevolumeyear)
\startpage--\finishpage\nl
Published: \publishdate}

\vglue .25truein

{\parskip=0pt\leftskip 0pt plus
1fil\def\\{\par\smallskip}{\Large\bf\thetitle}\par\medskip} \vglue
0.05truein

{\parskip=0pt\leftskip 0pt plus 1fil\def\\{\par}{\sc\theauthors}
\par\medskip}%
 
\vglue 0.03truein 

{\small\leftskip 25truept\rightskip 25truept{\bf Abstract}\stdspace\theabstract

{\bf AMS Classification}\stdspace\theprimaryclass
\ifx\thesecondaryclass\relax\else; \thesecondaryclass\fi\par
{\bf Keywords}\stdspace \thekeywords\par}\vglue 7truept

}   %%%% end of definition of \makeagttitle

\ifplaintex
\font\phead=cmsl9 scaled 950
\font\pnum=cmbx10 scaled 913
\font\pfoot=cmsl9 scaled 950
\headline{\vbox to 0pt{\vskip -4.5mm\line{\small\phead\ifnum
\count0=\startpage ISSN 1472-2739 (on-line) 1472-2747 (printed)
\hfill {\pnum\folio}\else\ifodd\count0\def\\{ }% 
\ifx\theshorttitle\relax\thetitle\else\theshorttitle\fi\hfill{\pnum\folio}
\else\def\\{ and }{\pnum\folio}\hfill\ifx\theshortauthors\relax\theauthors
\else\theshortauthors\fi\fi\fi}\vss}}
\footline{\vbox to 0pt{\vglue 0mm\line{\small\pfoot\ifnum\count0=\startpage
\copyright\ \gtp\hfill\else
\agt, Volume \thevolumenumber\ (\thevolumeyear)\hfill\fi}\vss}}
\else
\font=cmsl9 scaled 1050
\font\lnum=cmbx10 
\font=cmsl9 scaled 1050
\makeatletter
\def\@oddhead{{\small\lhead\ifnum\count0=\startpage ISSN 1472-2739 
(on-line) 1472-2747 (printed)\hfill {\lnum\number\count0}\else\ifodd\count0
\def\\{ }\ifx\theshorttitle\relax \thetitle \else\theshorttitle\fi\hfill
{\lnum\number\count0}\else\def\\{ and }{\lnum\number\count0}
\hfill\ifx\theshortauthors\relax 
\theauthors\else\theshortauthors\fi\fi\fi}}\def\@evenhead{\@oddhead}
\def\@oddfoot{\small\lfoot\ifnum\count0=\startpage\copyright\ \gtp\hfill\else
\agt, Volume \thevolumenumber\ (\thevolumeyear)\hfill\fi}
\def\@evenfoot{\@oddfoot}
\makeatother
\fi
\let\maketitlepage\makeagttitle
\let\makeshorttitle\maketitlepage
\let\maketitle\maketitlepage

\newwrite\gtoutfile
\long\gdef\makeheadfile{  %%% start of definition of \makeheadfile
{\def\\{, }\def\s{ }
\immediate\openout\gtoutfile head.xxx
\immediate\write\gtoutfile{Proxy-for: \ifx\theasciiauthors\relax
\theauthors\else\theasciiauthors\fi\s<\ifx\theasciiemail\relax\theemail\else\theasciiemail\fi>}
\immediate\write\gtoutfile{\noexpand\\}
\immediate\write\gtoutfile{Authors: \ifx\theasciiauthors\relax
\theauthors\else\theasciiauthors\fi}
{\def\\{ }\immediate\write\gtoutfile{Title: \ifx\theasciititle\relax
\thetitle\else\theasciititle\fi}}
\immediate\write\gtoutfile{Subj-class: GT or SG, GR etc}
\immediate\write\gtoutfile{MSC-class: \theprimaryclass\ifx\thesecondaryclass\relax\else, \thesecondaryclass\fi}
\immediate\write\gtoutfile{Journal-ref: Algebr. Geom. Topol. \thevolumenumber\s
(\thevolumeyear) \startpage-\finishpage}
\immediate\write\gtoutfile{Comments: Published by Algebraic and
Geometric Topology at}
\immediate\write\gtoutfile{\s\s\s  http://www.maths.warwick.ac.uk/agt/AGTVol\thevolumenumber/agt-\thevolumenumber-\thepapernumber.abs.html}
\immediate\write\gtoutfile{\noexpand\\}
\immediate\write\gtoutfile{}
\ifx\theasciiabstract\relax
\immediate\write\gtoutfile{\theabstract}\else
\immediate\write\gtoutfile{\theasciiabstract}\fi
\immediate\write\gtoutfile{}
\immediate\write\gtoutfile{\noexpand\\}
\immediate\write\gtoutfile{}
\immediate\closeout\gtoutfile}}  %%% end of definition of \makeheadfile

\def\maketitlepage{\makeagttitle\makeheadfile}
\let\makeshorttitle\maketitlepage
\let\maketitle\maketitlepage

%% file: agt-3-36.bbl
\begin{thebibliography}{99}
\bibitem[AK]{AK} S. Akbulut, H. King, {\sl Topology of Real Algebraic
    Sets.} L'Enseignement Math\'ematique, 29 (1983), 221-261.
\bibitem[Ba]{Ba} H. H. Baues, {\sl Obstruction Theory, on Homotopy
    Classification of Maps.} Springer Lectures Notes 628, 1977.
\bibitem[BR]{B-R} R. Benedetti and Jean-Jacques Risler, {\sl Real
    Algebraic and Semi-Algebraic Sets.} Hermann, \'Editeurs des
  science et des arts, 1990.
\bibitem[BT]{BT} T. Br\"ocker and  T. tom Dieck, {\sl Representations of
    Compact Lie Groups.} Springer-Verlag, 1985.
\bibitem[CF]{C-F} P. Conner and E. E. Floyd, {\sl Differentiable
    Periodic Maps.} Springer-Verlag, 1972.
\bibitem[F]{F} M. H. Freedman, {\sl The Topology of Four-Dimensional
    Manifolds.} J. Differ. Geom. 17 (1982), 357-453.
\bibitem[H]{H} M. W. Hirsch, {\sl Differential Topology.}
Springer-Verlag, 1976.
\bibitem[Hae]{Hae} A. Haefliger, {\sl Differentiable Imbeddings.}
Bull. Amer. Soc. 67 (1961), \mbox{109-112}.
\bibitem[Hi]{Hi} H. Hironaka, {\sl Resolutions of Singularities of an
    Algebraic Variety over a Field of Characteristic Zero: I,II.}
Ann. of Math. 79 (1964), 109 - 326.
\bibitem[Ke]{Ke} M. A. Kervaire, {\sl Some Non-Stable Homotopy Groups of
  Lie Groups.} Illinois J. Math. Vol 4 (1960), 161-169.
\bibitem[Kr1]{Kr1} M. Kreck, {\sl Surgery and Duality.} Annals of
  Mathematics Vol. 149 (1999), \mbox{704-754}.
\bibitem[Kr2]{Kr2} M. Kreck, {\sl Differential Algebraic Topology.}
Preprint.
\bibitem[Kr3]{Kr3} M. Kreck, {\sl Stratifolds.}
Preprint.
\bibitem[Mi1]{Mi_s} J. Milnor, {\sl A Procedure for Killing the
    Homotopy Groups of Differentiable Manifolds.} Symposia in Pure
  Math. A.M.S., Vol. III (1961), 39-55.
\bibitem[Mi2]{Mi0} J. Milnor, {\sl Morse Theory.} Annals of
  Math. Studies 51, Princeton University Press, 1965.
\bibitem[Mi3]{Mi_sp} J. Milnor, {\sl Remarks Concerning Spin Manifolds.}
Differ. and Combinat. Topology, Sympos.
Marston Morse, Princeton, 1965, 55-62.
\bibitem[Mi4]{Mi} J. Milnor, {\sl Singular Points of Complex
    Hypersurfaces.}
Princeton University Press, 1968.
\bibitem[St]{St} R. E. Stong, {\sl Notes on Cobordism Theory.} Mathematical Notes,
Princeton University Press, 1968.
\bibitem[Ste]{Ste} N.Steenrod, {\sl The Topology of Fibre Bundles.}
Princeton University Press, 1951.
\bibitem[Th]{Th} R. Thom, {\sl Quelques Properi\'et\'es Globales des
    Variet\'et\'es Differentiables.} Comm. Mathe. Helv. Vol 28 (1954),
  17-86.
\bibitem[W1]{W0} C. T. C. Wall, {\sl Classification of
    $(n-1)$-Connected $2n$-Manifolds.} Annals of Mathematics Vol. 75 (1962),  163-189.
\bibitem[W2]{W1}  C. T. C. Wall, {\sl Surgery on Compact Manifolds.}
  London Math. Soc. Monographs 1, Academic Press, New York, 1970.

\end{thebibliography}
